\documentclass[sort&compress]{elsarticle}

\usepackage{defns} 

\usepackage{pdflscape} 
\usepackage{afterpage}

\usepackage[utf8]{inputenc}
\usepackage{amsmath, amssymb,mathtools}
\usepackage{algorithm}
\usepackage{algpseudocode}
\usepackage{graphicx,appendix}
\usepackage{bm}
\usepackage{tikz}

\usetikzlibrary{shapes,arrows}

\usepackage{wrapfig}

\usepackage{amsthm}
\newtheorem{example}{Example}

\newcommand{\R}{\mathbb{R}}
\newcommand{\M}{\mathbb{M}}

\newcommand{\Prob}{\mathbb{P}}

\newcommand{\T}{\mathcal{T}}

\newcommand{\tr}{\text{tr}}
\newcommand{\cmt}[1]{} 

\title{A Metric Tensor Approach to Data Assimilation with Adaptive Moving Meshes
\tnoteref{label1}}
\tnotetext[label1]{This research was supported in part by NSF
grants DMS-1714195 and DMS-1722578.}
\author[1]{Cassidy Krause}
\ead{ckrause@ku.edu}

\author[1]{Weizhang Huang} \ead{whuang@ku.edu}

\author[2]{David B Mechem} \ead{dmechem@ku.edu}

\author[1]{Erik S Van Vleck}\ead{erikvv@ku.edu}

\author[3]{Min Zhang}\ead{minzhang@math.pku.edu.cn}

\address[1]{Department of Mathematics, University of Kansas, Lawrence, KS
}
\address[2]{Department of Geography and Atmospheric Science, University of Kansas, Lawrence, KS}

\address[3]{School of Mathematical Sciences, Peking University, Beijing 100871, China}

\date{\today}

\tikzstyle{block} = [rectangle, draw, fill=blue!20,
    text width=4cm, text centered, rounded corners, minimum height=4em]
\tikzstyle{line} = [draw, -latex']
\tikzstyle{cloud} = [draw, ellipse,fill=red!20, text width = 2.5cm, text centered, node distance=4cm,
    minimum height=4em]

\begin{document}

\begin{abstract}
Adaptive moving spatial meshes are useful for solving physical models given by time-dependent partial differential equations. However, special consideration must be given when combining adaptive meshing procedures with ensemble-based data assimilation (DA) techniques. In particular, we focus on the case where each ensemble member evolves independently upon its own mesh and is interpolated to a common mesh for the DA update. This paper outlines a framework to develop time-dependent reference meshes using locations of observations and the metric tensors (MTs) or monitor functions that define the spatial meshes of the ensemble members. We develop a time-dependent spatial localization scheme based on the metric tensor (MT localization). We also explore how adaptive moving mesh techniques can control and inform the placement of mesh points to concentrate near the
location of observations, reducing the error of observation interpolation. This is especially beneficial when we have observations in locations that would otherwise have a sparse spatial discretization. We illustrate the utility of our results using discontinuous Galerkin (DG) approximations of 1D and 2D inviscid Burgers equations. 
The numerical results show that the MT localization scheme compares favorably with standard Gaspari-Cohn localization techniques. In problems where the observations are sparse, the choice of common mesh has a direct impact on DA performance. The numerical results also demonstrate the advantage of DG-based interpolation over linear interpolation for the 2D inviscid Burgers equation.
\end{abstract}

\maketitle

\section{Introduction}\label{sec:Intro}

Data assimilation (DA) combines noisy data, usually from instrumental uncertainty or issues with scale, 
with models that are imperfect, due to simplifying assumptions or other inaccuracies, to improve predictions about the past, current, or future state of a system. Typical DA techniques require a prior uncertainty and use the likelihood of observations to calculate the posterior distribution. This Bayesian context provides not only predictions but also quantification of the uncertainty in these predictions. 
DA originated with numerical weather prediction and is now employed in many scientific and engineering disciplines. With infinite dimensional models such as partial differential equations (PDEs) there are often features that require fine spatial resolution to obtain accurate approximations of model solutions. An alternative to uniform fine meshes is to employ meshes that are coarse in parts of the spatial domain but fine in other parts of the domain. When the regions requiring finer resolution change with time, time-dependent adaptive meshes are advantageous. In the context of data assimilation, employing time-dependent adaptive meshes also has the potential to reduce the dimension needed for good approximation of model solutions.

Combining DA with adaptive moving mesh techniques presents both opportunities and challenges. Many DA techniques employ an ensemble of solutions to make predictions and estimate the uncertainty in the predictions.
This introduces a key choice for implementation: either each ensemble solution may evolve on its own independent adaptive mesh, or the ensemble time evolution can happen on a mesh that is common to all ensemble members. 
In the first case, the mesh for each ensemble member may be tailored to each ensemble solution. However, when incorporating a DA scheme, this presents the challenge of combining ensemble solutions that are supported on different spatial meshes.
In the second case, one needs to combine the features of the different ensemble solutions to determine a mesh that provides, on average, good approximation for all ensemble members. In addition, the relative position of the ensemble mesh(es) with respect to potentially time-dependent observation locations may also motivate the choice of ensemble mesh(es).


Our contribution in this paper is to develop a framework and techniques to utilize adaptive meshing techniques for data assimilation with application to PDE models in one and higher space dimensions. This framework allows each ensemble member to evolve on its own independent mesh and presents an adaptive common mesh for the DA update.
An adaptive mesh, while not usually uniform in the standard Euclidean metric, 
can be viewed as uniform in some other metric. This metric is defined by a positive-definite matrix valued monitor function, also called a metric tensor, which controls mesh movement for the ensemble members and is also used to define the new adaptive common mesh. The computation of the metric tensor easily generalizes to higher spatial dimensions. Using an adaptive common mesh based on metric tensors of the ensemble meshes provides a means for determining a mesh that is common to all ensemble members while providing good approximation properties for the ensemble solutions.

The adaptive common mesh is formed through the metric tensor intersection of the metric tensors that monitor the movement of the ensemble meshes. Geometrically, the metric tensor intersection is the same as circumscribing an ellipsoid on an element of two meshes and then finding an ellipsoid that resides in the geometric intersection of the first two ellipsoids. The new element given by that intersecting ellipsoid is the result of the metric tensor intersection. This procedure must be done pairwise, so the resulting ellipsoid is not necessarily maximal. However, a greedy algorithm can be used to find an ordering that seeks to maximize the resulting ellipsoid. 



The metric tensor intersection of the ensemble member meshes forms a common mesh that supports all ensemble members with high accuracy. However, the mesh points of this common mesh may not align with the observation locations. If the observation locations do not coincide with nodes in the common mesh, then the location mismatch results in errors from interpolating the observations to the common mesh. This is especially relevant in the case of time-dependent observation locations. 
With field alignment \cite{Emanuel2007}, a variational DA scheme is developed that assimilates both the alignment and amplitude of observations simultaneously. In addition, a two-step approach is developed, where first the locations of the state variables are adjusted to better match the observations. Rather than assume that the observations occur at the same spatial discretization used for the numerical solution, a vector of displacement is employed so that the (interpolated) numerical solution at adjusted nodes is obtained to maximize the posterior distribution of the numerical solution with displacement given the observations. The DA then proceeds with traditional correction based upon the amplitude of the errors of the numerical solution at the adjusted discretization.

The mismatch of observation locations and nodal positions of the common mesh presents another potential opportunity for DA on adaptive moving meshes: the meshes can adapt to concentrate near or align with the observation locations. 
An observation mesh can be formed by associating a metric tensor with the location of the (potentially time-dependent) observations. Intersecting the common mesh from the ensemble members with this observation mesh provides a new common mesh that is concentrated near observation locations. Interpolation error can have a significant effect on the accuracy of DA schemes, and concentrating the mesh near observation locations reduces the interpolation required, thereby improving the performance of the DA algorithms. Of course, a fixed common mesh can also be used to concentrate the mesh near fixed observation locations, but this approach has the benefit of adapting easily to time-dependent observation locations.

In addition, we develop spatially and temporally dynamic localization schemes, based upon the metric tensor(s) corresponding to the adaptive common mesh. Localization improves DA procedures by ensuring that observations only affect nearby points. Broadly speaking, localization schemes fall into two categories: domain localization and covariance (or R) localization. Domain localization schemes define a spatial radius and use that to define which mesh points are affected by a given observation. Covariance localization schemes use a correlation function to modify the covariance matrix that is used in the DA update, so that the covariance between an observation and the solution values decays to zero as the distance between the observation and the solution values increases.

We develop a domain localization scheme that employs the metric tensor. Employing a fixed, uniform radius of influence for the observations is that the localization scheme may not be effective if there is a steep gradient in the solution. One could predetermine the location of the gradient and adjust the localization scheme accordingly, but if the regions of large gradient are time-dependent, this will usually result in the tuned localization parameter being quite small. However, since the metric tensor provides information about the dynamics of the ensemble solution, it can be used to define an adaptive localization scheme where the localization radius can vary in time and space.

One benefit of using an adaptive moving mesh is that fewer mesh points can be used while still maintaining the same accuracy. Having an adaptive time-dependent common mesh allows for a fewer number of nodes used in the common mesh as compared to a fixed, fine common mesh, increasing the efficiency of the linear algebra, e.g., when updating the mean and covariance with an ensemble Kalman filter. An efficient implementation of the Ensemble Kalman Filter (EnKF) requires $\mathcal{O}((D+N_e)D^2 + (M+D)N_e^2)$ flops when $D\ll M$ or more generally, for example when $D\approx M$, $\mathcal{O}((M+D+N_e)N_e^2)$ (see, e.g., \cite{Mandel2006}) where $D$ is the dimension of the observation space, $M$ is the dimension of the discretized dynamical system, and $N_e$ is the number of ensemble members.
In large scale geophysical applications we typically desire $N_e\approx 20$ (in general $N_e$ should be roughly the number of positive and neutral Lyapunov exponents).
A reduction in $M$ based upon using fewer mesh points while maintaining or enhancing accuracy results in improved efficiency.

There are several recent works on integrating adaptive spatial meshing techniques with DA, although most of the focus has been on PDE models in one space dimension. This includes methods based on evolving meshes based on the solution of a differential equation, methods in which meshes are updated statically based upon interpolation, and remeshing techniques that add or subtract mesh points as the solution structure changes.
In \cite{Bonan}, the evolution of meshes was movement of the nodes was determined by the solution of moving mesh differential equations that are coupled to the discretized PDE. The state variables of the PDE were then augmented with the position of the nodes and incorporated into a DA scheme. The test problem consisted of a two-dimensional ice sheet assumed to be radially symmetric; therefore, it reduced to a problem with one spatial dimension.
In \cite{Colin} and \cite{Du} common meshes were developed based on a combining through interpolation the ensemble meshes. This allowed for update of mean and covariance for Kalman filter based DA techniques while allowing each ensemble member to evolve on its own independent mesh. That is, at each observational timestep, the ensemble members were interpolated to the common mesh, updated with the DA analysis, and then interpolated back to their respective meshes. Specifically, a uniform, non-conservative mesh was used in \cite{Colin}, with Lagrangian observations in one spatial dimension. Higher spatial dimensions were used in \cite{Du}, with a fixed common mesh refined near observation locations. \cite{sampson2020} uses the same 1D non-conservative adaptive meshing scheme as in \cite{Colin} and extends this approach through the use of an adaptive common mesh, where, like in \cite{Bonan}, the state vector is augmented with the node 
locations.

The outline of this paper is as follows. Background of data
assimilation and adaptive moving mesh techniques is given in Section
\ref{sec:background}. This includes the framework we develop to include equations describing mesh movement within a DA framework. The development of adaptive meshing techniques for DA is
detailed in Section \ref{sec:adaptiveDA}. Metric tensors are introduced and their connection to non-uniform meshes is discussed. Techniques for combining meshes based on metric tensor intersection and for concentrating ensemble mesh(es) near observation locations are developed. 
The details of our implementation are in Section \ref{sec:details}. This includes the discontinuous Galerkin discretization we employ and the specific metric tensor formulation we use to adaptively evolve the ensemble mesh(es). The details of our experimental setup and numerical results for both 1D and 2D inviscid Burgers equations are presented
in Section \ref{sec:results}.

\section{Background on Data Assimilation and Adaptive Moving Meshes}\label{sec:background}

\subsection{Data Assimilation} Data assimilation techniques seek to combine models and data to improve predictions and quantify uncertainty typically in a Bayesian context (see, e.g., \cite{Carrassi2018}, \cite{LSZ}, \cite{ReichCotter15}, \cite{VanLeeuwen18}).
Consider a finite dimensional discrete time system for a state vector $u \in \R^M$ ($M > 0$) that evolves based upon
\begin{equation}\label{eq:dtime_mod}
   u_{n+1} = \Psi(u_n) + \xi_{n}, \quad \xi_n \sim \mathcal{N}(0,\Sigma),
\end{equation}
where $\Psi(\cdot)$ propagates the state forward in time, $n$ stands for the $n$th time step, and $\xi_n$ is assumed to be a normally distributed model error with covariance matrix $\Sigma$ and mean $0$. Equation \eqref{eq:dtime_mod} can be used to forecast the state of the dynamical system, but this prediction can often be improved by including the observation $y_{n+1} \in \R^D$ ($D > 0$)
given by the data model
\begin{equation} \label{eq:dtime_obs}
    y_{n+1} = \mathcal{H}(u_{n+1}^t)+\eta_{n+1}, \quad \eta_{n+1} \sim \mathcal{N}(0,R),
\end{equation} where $u_{n+1}^t$ is the unknown ``truth'' at time $n+1$, $\mathcal{H}: \R^M \rightarrow \R^D$ is the observation operator, and $\eta_{n+1}$ is assumed to be a normally distributed observation error with covariance matrix $R$ and mean $0$. In many applications, $D \ll M$. We wish to determine $\{u_n\}_{n=0}^N$ that satisfies, in some sense, both the physical and data models. Note that implicit in the formulation of the data model \eqref{eq:dtime_obs} is that the state variables ($u_{n+1}$) and data variables ($y_{n+1}$) are supported on common spatial locations.


Many data assimilation techniques are based upon a Bayesian approach that determine the posterior distribution from the prior distribution and the likelihood. Given an observation $y_n$ at time $t_n$ and a prior estimate
$\Prob(u_n)$ of the state, Bayes' Theorem states that
\begin{equation}\label{eq:bayes}
\Prob(u_n|y_n)\, { \propto }\, \Prob(y_n|u_n)\Prob(u_n).
\end{equation}
This procedure extends to the sequential assimilation of observations at multiple times under the assumption that the state is Markovian. Note that since the model noise $\{\xi_j\}_{j=1}^{N}$ is independent and identically distributed (i.i.d.), the prior can be written as $\Prob(u_{n}) = [\prod_{j=1}^{n}\Prob(u_{j}|u_{j-1})]\Prob(u_0)$, where, e.g., $u_0 \sim \mathcal{N}(u_0^b,P_0^b)$.



For nonlinear models \eqref{eq:dtime_mod} available data assimilation techniques include the ensemble Kalman filter (EnKF), particle filters (PF), variational methods such as 4DVar, and hybrid techniques that seek to combine the best features of different types of techniques. Several of these classes of methods are naturally ensemble based (EnKF, PF and their variants, as well as several hybrid methods) while variational methods may employ ensembles of solutions to approximate derivatives or as part of an iterative linear system solver.



Ensemble DA procedures use an ensemble of solutions to make predictions for the physical state via a two-step process. First, the prediction step uses the physical model \eqref{eq:dtime_mod} to integrate the ensemble members $\{u_n^{e_i}\}_{i=1}^{N_e}$, where $N_e$ is the number of ensemble members, to make the ensemble forecasts $\{\hat u_{n+1}^{e_i}\}_{i=1}^{N_e}$. Second, the analysis step incorporates the observations at time $t_{n+1}$ to adjust the prediction $\{u_{n+1}^{e_i}\}_{i=1}^{N_e}$. EnKF does this through the following sequential process:
\begin{equation}
    \text{Prediction:} \quad \begin{cases}\hat u_{n+1}^{e_i} &= \Psi(u_n^{e_i})+\xi_n^{e_i},\quad e_i = 1, \dots, N_e\\
    \hat m_{n+1} &= \frac{1}{N_e}\sum_{i=1}^{N_e} \hat u_{n+1}^{e_i}, \\
    P^b_{n+1} &= \frac{1}{N_e-1}\sum_{i=1}^{N_e}\left(\hat u_{n+1}^{e_i}-\hat m_{n+1}\right)\left(\hat u_{n+1}^{e_i}-\hat m_{n+1}\right)^T ,
    \end{cases} \label{eq:DA_pred}
\end{equation}
\begin{equation}
    \text{Analysis:} \quad \begin{cases} \mathbf{K}_{n+1} &= P^b_{n+1}H^T\left(HP^b_{n+1} H^T + R\right)^{-1}, \\
    u_{n+1}^{e_i} &= \left(I - \mathbf{K}_{n+1}H\right)\hat u_{n+1}^{e_i} + \mathbf{K}_{n+1}y_{n+1}, \quad e_i = 1, \dots, N_e
    \end{cases}\label{eq:DA_analysis}
\end{equation}
where $I$ is the identity matrix, $H$ is the linearization of $\mathcal{H}$, and $\mathbf{K}_{n+1}$ is the Kalman gain matrix. 
 The ensemble mean $\hat m_{n+1}$ together with the forecast covariance $P_{n+1}^b$ in \eqref{eq:DA_pred} is employed to make predictions with a measure of uncertainty. Again, the implicit assumption here is that all ensemble members reside on the same mesh. 

Variational methods such as 4DVar and 3DVar are direct formulations of Bayes' Theorem.
The analysis update of a variational DA algorithm can be viewed as the minimizer of a corresponding cost function. For example, when the prior distribution and the observational error model are Gaussian, the cost function is
\begin{equation}
\label{4DVarCF}
    J(u_0) = (u_0-u_b)^TB^{-1}(u_0-u_b) + \sum_{n=0}^N (y_n - \mathcal{H}(u_{n}))^T R_n^{-1}(y_n - \mathcal{H}(u_{n})),
\end{equation}
where $u_n$ is obtained by the evolution of the model dynamics \eqref{eq:dtime_mod} and $B$ is the background error covariance. 
Note that this cost function is quadratic when both the physical model $\Psi$ and the observation operator $\mathcal{H}$ are linear, so in this case there is a unique optimizer. In the case that either the physical model or observation operator are nonlinear, or if the error distributions are not Gaussian, the posterior distribution may not have a unique optimum. The use of non-Gaussian error distributions results to corresponding terms in the cost function \eqref{4DVarCF} which may be approximated at least locally with a quadratic cost function, e.g., with variants such as incremental 4DVar.

Both EnKF and 4DVar are used in large-scale applications. 4DVar provides great flexibility in including terms in the cost function, does not require linear physical model or observation operators, and can make use of sophisticated optimization techniques. 
EnKF is based upon Gaussian assumptions and parameterizes the prediction and uncertainty in terms of the ensemble mean and sample covariance. Together with local linear approximation of the observation operator, EnKF is a linear solver corresponding to a quadratic cost function.
For more on the advantages and disadvantages of the EnKF and 4DVar techniques, see \cite{Kalnay2007, Lorenc2003}.

Hybrid methods seek to combine the advantages of ensemble Kalman filter techniques and variational techniques.
An important motivation behind hybrid methods
is to incorporate flow-dependent background error covariance matrices ($P_{n+1}^b$) into a variational setting.
A representative hybrid method is the ETKF-4DVar technique in which $B$ in \eqref{4DVarCF} is replaced by
\[
\tilde B = \beta B + (1-\beta)P_b
\]
where $0\leq \beta \leq 1$ is a weighting factor, $B$ is the background error covariance used in 4DVar, and $P_b$ is the error covariance found using the ensemble transform Kalman filter (ETKF), a square root filter similar to the EnKF \eqref{eq:DA_pred}, \eqref{eq:DA_analysis}. For $\beta =1$ this reduces to the standard 4DVar, while for $\beta=0$ it is known as
4DVAR-BEN, and for $\beta=1/2$, the method is known as the ETKF-4DVAR. For further discussion of hybrid methods, see, e.g., \cite{Asch2016, Bannister2017, Lorenc2013}.

While the previous classes of techniques all rely to some extent on Gaussian assumptions by parameterizing the predicted state and its uncertainty in terms of a mean and covariance, particle filters approximate the posterior distributions in an unstructured way in terms of particles (analogous to ensemble members) and particle weights. For example, the standard or bootstrap particle filter uses the model dynamics \eqref{eq:dtime_mod} to make predictions using each particle.
The weights $\{ w_{n-1}^{e_i}\}_{i=1}^{N_e}$ are updated using the observation $y_n$. In the bootstrap PF the update is
\begin{equation}\label{filtStep}
w_n^{e_i} = c\, w_{n-1}^{e_i} \Prob(y_n| u^{e_i}_n),
\end{equation}
where the likelihood, given Gaussian observational error model, is
\begin{align}
\label{likelihood}
\Prob(y_n|u_n) \propto \exp\left[-\frac{1}{2} \left(y_n - H u_n\right)^T R^{-1} \left(y_n - H u_n\right) \right]
\end{align}
and $c$ is chosen so that $\sum_{i=1}^{N_e} w_{n}^{e_i}=1$. The standard particle filter and many variants suffer from weight degeneracy in which most weights are nearly zero and this decreases the effective number of particles. This leads to the need for a large number of particles ($N_e$) that increases exponentially with the dimension of the state space $(M)$ and the observation space $(D)$. Several variants including the implicit PF, the equivalent weights PF, and optimal proposal PF have been developed (see, e.g., \cite{VanLeeuwen18}) to overcome weight degeneracy and the curse of dimensionality in PFs.


Modern DA procedures employ spatial localization and covariance inflation techniques. Localization can be done either through a Schur product (also known as a Hadamard or element-wise product) with the covariance matrix, resulting in observations having little or no effect on distant nodes. It can also be done through domain localization, where the domain is broken into subdomains, and the observations of each subdomain only affects the variables supported within that subdomain. Inflation (either additive or multiplicative) increases the entries of the covariance matrix, preventing degeneracy of the procedure.

The focus of this work is employing adaptive moving meshes with ensemble DA methods. In the following we outline a framework for reducing to the discrete time, finite dimensional DA problem \eqref{eq:dtime_mod}, \eqref{eq:dtime_obs}. We start from the simplest cases of discretized ODE and PDE on fixed spatial meshes,
and then formulate the DA problem using moving mesh methods in which a differential equation determines the movement of mesh points. We will emphasize the errors that are introduced in these different scenarios.

\subsection{ODEs and Time-Dependent PDEs Discretized on Fixed Mesh}


Given an ODE
\begin{equation}
\frac{du}{dt} = F(u,t)
\label{eq:mod_ode}
\end{equation}
as the physical model, the use of a time stepping technique and the subsequent addition of additive model error is used to obtain \eqref{eq:dtime_mod}. Similarly, given a time-dependent PDE, after spatial discretization we obtain a system of ODEs that can then be discretized in time to obtain a discrete time, finite dimensional physical model \eqref{eq:dtime_mod} with the addition of model noise.

\subsection{The Physical PDE Model and Moving Mesh Equations}
Consider now the physical model as a time-dependent PDE written abstractly as $\frac{\partial u}{\partial t} = F(u)$
posed on an appropriate function space.
PDEs employed in a DA context may be coupled to an equation that evolves the spatial mesh, enabling DA on an adaptive moving mesh. There are two basic approaches to adaptive moving meshes. The first is a rezoning approach, used in \cite{Colin,sampson2020}, which updates the mesh at each time using a given mesh generation and interpolates the PDE solution from the old mesh to the new mesh. The second approach is a quasi-Lagrange approach where the mesh is considered to move continuously with time. In this case, the discretized PDE is supplemented with an advective term to reflect mesh movement.
If $u = u(x(t),t)$ satisfies a PDE given by $u_t = F(u)$, then the total derivative $\dfrac{du}{dt}$ is given by
\begin{equation}
    \frac{du}{dt} = \nabla_x u \cdot \frac{dx}{dt} + u_t =  \nabla_x u \cdot \frac{dx}{dt} + F(u) \equiv {\mathcal F}(u,x) .\label{eq:mmpde}
\end{equation}

The equation for the mesh movement comes from a variational approach in which a cost function (the meshing functional) is minimized
through a gradient flow differential equation:
\begin{equation}
    x_t = - \frac{1}{\tau}\nabla_x \mathcal{L}(x,u)\\
    \equiv {\mathcal G}(x,u), \label{eq:mesheq}
\end{equation}
where $\tau\geq 0$ is a user-specified parameter controlling the speed of mesh movement. We note that when $\tau=0$, \eqref{eq:mesheq} reduces to an algebraic equation ${\mathcal G}(x,u) = 0$ and when satisfied the mesh instantaneously satisfies a local minimizer of the meshing functional.
In practice,
first fix $x$ and then update the solution of the PDE by solving \eqref{eq:mmpde}. The solution to \eqref{eq:mmpde}, together with $x$, gives the mesh velocity defined by \eqref{eq:mesheq}.
A new mesh is obtained by integrating \eqref{eq:mesheq} in time.



In an ensemble-based method, an adaptive mesh PDE discretization can be used for each of the ensemble members. However, the computations of the ensemble mean and covariance only make sense if the values of the ensemble members are taken from the same spatial locations. Therefore, if the ensemble members' meshes can evolve independently via an adaptive moving mesh scheme, special care must be taken to calculate the mean and covariance at each DA step. Previous works have explored two general approaches to this problem. One approach is to interpolate the ensemble solutions to a common mesh at each observational time step and assimilate the PDE variables only. The common mesh approach is used here, as well as in \cite{Colin,Du}. Another method is to assimilate both the PDE variables and the common mesh locations, as done in \cite{Bonan}.

In \cite{Bonan}, they augment the state vector $U^n$
with the mesh points $x^n$. That is, the new state vector is given by
\begin{equation}
    V^n = \begin{bmatrix} U^n \\ x^n \end{bmatrix}.
\end{equation}
The DA update is then performed on this augmented state vector. Note that this approach doubles $M$ in the DA scheme,
affecting the computational cost of the DA update.
In \cite{Colin,sampson2020} the mesh movement is based on a Lagrangian flow of form $\frac{dx}{dt} = u$
together with upper and lower bounds on the mesh spacing that are enforced by remeshing and subtracting
or adding mesh points.
For example, on a spatial domain $\Omega = (0,1)$ this corresponds
to a meshing functional \eqref{eq:mesheq} of the form
\begin{equation}
     \mathcal{L}(x,u) = -\int_0^x u(y,t)dy.
\end{equation}

\begin{example} As a concrete example of a PDE and mesh movement equation,
consider a 1D reaction diffusion equation:
\begin{equation}
    u_t - u_{xx} = f(u),\quad 0<x<1,\quad t>0
\end{equation}
with, for example, homogeneous Dirichlet or homogeneous Neumann boundary conditions and smooth initial data.

A non-uniform finite difference discretization in space yields $(U_j\equiv U_j(t) \approx u(x_j,t))$:
\begin{equation}
\label{eq:int:grid:update:mmpde5:eqU}
\begin{array}{lcl}
\dot{U}_j & = & \Big[ \frac{U_{j+1} - U_{j-1}}{x_{j+1} - x_{j-1}} \Big] \dot{x}_j +    \frac{2}{x_{j+1} - x_{j-1}}
\big[ \frac{U_{j-1} - U_j}{x_{j} - x_{j-1} } + \frac{U_{j+1} - U_j}{x_{j+1} - x_j} \big]
+ f( U_j).
\end{array}
\end{equation}
Now consider a moving mesh scheme determined by the equidistribution of arc length of the solution:
\begin{equation}
\label{eq:int:grid:update:mmpde5:eqX}
\tau \dot{x}_j  =
\sqrt{  (x_{j+1} - x_j)^2 + (U_{j+1} - U_j)^2 }
- \sqrt{ (x_{j-1} - x_j)^2 + (U_{j-1} - U_j)^2  } .
\end{equation}
After discretizing both
 \eqref{eq:int:grid:update:mmpde5:eqU} and \eqref{eq:int:grid:update:mmpde5:eqX} in time we obtain a discrete time, finite dimensional model \eqref{eq:dtime_mod} with the state vector $ U$ or $\begin{bmatrix} U \\ x \end{bmatrix}$ depending on whether only the PDE solution is assimilated or if the mesh locations are as well.
\end{example}

\subsection{Observations and Their Locations}
Consider $D$ observations at a fixed time $t$, $\mathbf{y}(t) = (y_1(t), ..., y_D(t))^T \in \R^D$ supported on a potentially time-dependent set of observation locations $\{x_j^o(t)\}_{j=1}^{D}$. In practice, these observations are related to the unknown truth $\mathbf{u^t}$ through the observation operator, suppressing the time dependence,
\begin{equation} \label{eq:obs_time_dependent}
\mathbf{y} = \mathcal{H}(\mathbf{u^t}) + \eta, \qquad \eta\sim\mathcal{N}(0,R).
\end{equation}
At a fixed time $t = t_n$, if $\{x_j^o(t_n)\}$ are a subset of the nodes of the common mesh, then \eqref{eq:obs_time_dependent} is exactly \eqref{eq:dtime_obs}, and the DA update proceeds as detailed above.
However, it is often the case that the $\{x_j^o(t_n)\}$ are not a subset of the nodes of the common mesh. In this case, it is natural to either interpolate the ensemble members onto a common mesh that contains the observations or interpolate the observations to the common mesh to obtain innovations of the form
\begin{equation}
\mathbf{y} = \mathcal{H}(\mathcal{I}_{co}(\mathbf{u^{e_i}})) + \eta\,,\quad {\rm or }\quad
\mathcal{I}_{oc}(\mathbf{y}) = \mathcal{I}_{oc}[\mathcal{H}(\mathcal{I}_{co}(\mathbf{u^{e_i}})) + \eta] ,
\end{equation}
where $\mathcal{I}_{co}$ interpolates from the common mesh to the locations in the domain of the observation operator and $\mathcal{I}_{oc}$ maps from the observation locations to (part of) the common mesh.


\section{Development of Adaptive Mesh Methods for DA}\label{sec:adaptiveDA}

\begin{figure}
\begin{center}
\begin{tikzpicture}[node distance=2cm, auto]
\node [block] (M_int) {Metric tensor intersection of ensemble meshes};
\node [block, right of=M_int, node distance=4.5cm] (concen) {Concentrate mesh near observations};
\node [block, right of=M_int, below of=M_int, node distance = 3cm] (common) {Common mesh};
\node [block, below of=common] (interp1) {Interpolate ensemble members to common mesh};
\node [cloud, left of=interp1](loc) {Adaptive localization};
\node [block, below of=interp1] (DA) {DA update on common mesh};
\node [block, below of=DA] (interp2) {Interpolate ensemble members back to ensemble meshes};
\node [cloud, left of=interp2] (remesh) {Remesh ensemble};
\node [block, below of=interp2] (step) {Step forward ensemble members and meshes};
\path [line] (M_int) -- (common);
\path [line] (concen) -- (common);
\path [line] (common) -- (interp1);
\path [line,dashed] (common) -| (loc);
\path [line,dashed] (loc) |- (DA);
\path [line] (interp1) -- (DA);
\path [line] (DA) -- (interp2);
\path [line,dashed] (interp2) -- (remesh);
\path [line,dashed] (remesh) |- (step);
\path [line] (interp2) -- (step);
\end{tikzpicture}
\end{center}
\caption{Algorithm for DA on an adaptive moving mesh. Optional steps are in red; required steps are in blue.}
\label{fig:algorithm}
\end{figure}
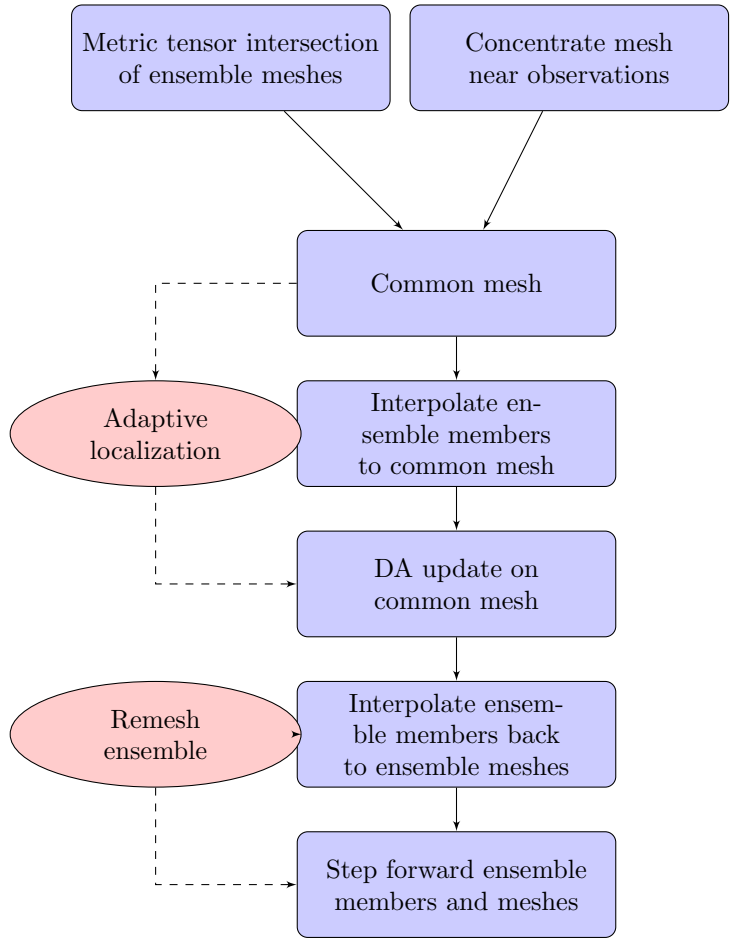

In this section we develop techniques for DA with adaptive moving meshes. These techniques are based upon the use of a metric tensor that describes the mesh. In particular, a non-uniform mesh is uniform with respect to the metric tensor being employed. The use of metric tensors is applicable not only in one spatial dimension, but also in higher spatial dimensions. We first describe the development of metric tensors using the so-called equidistribution and alignment conditions. Metric tensor intersections are introduced next and are used to combine meshes (e.g., the meshes of ensemble members) and potentially the locations of observations into what is in some sense an ``averaged'' mesh. We develop techniques for defining metric tensors that will concentrate mesh points near observation locations or other locations of interest and develop adaptive localization techniques based upon the use of metric tensors. A rough outline of the overall algorithm in given in Figure \ref{fig:algorithm}.


\subsection{Metric Tensors and Adaptive Meshes} \label{sec:m}
In one dimension, a mesh can be defined by determining the size of each of the elements of a mesh. In higher spatial dimensions, however, it is necessary to have conditions to also determine the shape and orientation of each element in addition to its size. For this and following discussions of adaptive mesh techniques, we follow the approach of \cite{Huang2010}.

Consider a polyhedral domain $\Omega \subset \R^d$ ($d\geq 1$), and let $\hat K$ be a reference element such that it is an equilateral $d$-simplex with unit volume. Then, given a simplicial mesh $\T_h$, for any element $K \in \T_h$, there is a unique invertible affine mapping $F_K: \hat K \rightarrow K$ such that $K = F_K(\hat K)$. The size, shape, and orientation of $K \in \T_h$ can be obtained from $F_K'$.

One way to create adaptive meshes is through the use of a given symmetric positive-definite matrix-valued monitor function, $\M = \M(x)$, which is also sometimes called a metric tensor. This monitor function defines a metric, and a mesh that is uniform in this metric is said to be an $\M$-uniform mesh. That is, a mesh $\T_h$ is $\M$-uniform if and only if all elements have a constant volume and are equilateral in the metric $\M$. These conditions are called the equidistribution and alignment criteria (\cite{Huang2010}, section 4.1.1). To be specific, let $\M_K$ to be the average of $\M(x)$ over the element $K$, that is,
\begin{equation}
    \M_K = \frac{1}{|K|}\int_K\M(x)dx.
\end{equation}
Then the equidistribution condition is given by
\begin{equation} \label{eq:equidistribution}
    \sqrt{\det (\M_K)}\; |K| = \frac{\sigma_h}{N}, \quad \forall K \in \T_h,
\end{equation}
where $\sigma_h = \sum_{K\in \T_h}|K|\sqrt{\det(\M_K)}$ and $N$ is the number of simplexes/elements in $\T_h$. The alignment condition is equivalent to
\begin{equation}\label{eq:alignment}
\frac{1}{d}\tr\left(\left(F_K'\right)^{-1}\M_k^{-1}\left(F_K'\right)^{-T}\right)
=\det\left(\left(F_K'\right)^{-1}\M_K^{-1}\left(F_K'\right)^{-T}\right)^{\frac{1}{d}}, \quad \forall K\in\T_h
\end{equation}
where $\det(\cdot)$ and $\tr (\cdot)$ denote the determinant and trace of a matrix, respectively.
Given a user-prescribed monitor function $\M = \M(x)$, equations \eqref{eq:equidistribution} and \eqref{eq:alignment} can be used to define an energy functional, which in turn generates mesh movement; see detail in Section~\ref{sec:details-mm}.

\subsection{Forming the Common Mesh} \label{sec:commonmesh}
At each data assimilation step, the common mesh is calculated by taking the
metric tensor intersection of the monitor functions for each ensemble
member; see, e.g., \cite{ZhangChengHuangQiu2018}.
The metric tensor intersection is better understood in geometry as the intersection of unit balls in different norms that are associated with symmetric and positive definite (SPD) matrices.
Mathematically, if $A$ and $B$ are SPD matrices of the same dimension,
then their intersection is denoted by ``$\cap$'' and defined as
\begin{equation}\label{eq:intersection}
A\cap B = P^{-1} \text{diag}\left(\max(1, b_1), \cdots,\max(1, b_d)\right) P^{-T},
\end{equation}
where $P$ is a nonsingular matrix such that
\begin{equation}
P A P^T = I,\qquad\qquad
P B P^T = \text{diag}(b_1,\cdots, b_d).
\end{equation}
It is not difficult to show that the unit ball in the norm associated with $A\cap B$ is contained in the unit balls in the norms associated with $A$ and $B$. As a consequence, when $A$ and $B$ are two metric tensors (meaning they are SPD-matrix-valued functions), the mesh associated with $A\cap B$ will combine concentration patterns of the meshes associated with $A$ and $B$.

Consider a DA scheme with $N_e$ ensemble members $u^{e_i}$, $i=1,...,N_e$ which are defined on different, independent meshes $\T^{e_i}_h$, $i=1,...,N_e$. Assume that the ensemble meshes $\T^{e_i}_h$ and the common mesh $\T_h$ have the same number of mesh elements and vertices (denoted $N$ and $M$, respectively), as well as the same connectivity. That is, they differ only in the location of the vertices.
At each observational time step, first interpolate the ensemble members from the ensemble meshes $\T^{e_i}_h$ to the current common mesh $\T_h$. Then compute the metric tensor for each of the ensemble members, denoted $\M^{e_i}_K$, and obtain a single metric tensor $\M_K^m$ by matrix intersection,
%
$\M_K^m = \M^{e_1}_{K}\cap\M^{e_2}_{K}\cap...\cap\M^{e_{N_e}}_{K}$,
%
which will be used to generate the common mesh at the next time step.
In practice, this computation can be done sequentially. In 1D, the order does not matter since the metric tensor intersection corresponds to finding the maximum value.
However, the order does matter in multi-dimensions and different orderings lead to different final metric tensors. While the examples we present in Section \ref{sec:results} are robust with respect to the different orderings, we have used the Greedy Algorithm based on minimizing the determinant (which is equivalent to maximizing the area of the ellipsoid from the intersection) to determine an optimal ordering.

\subsection{Concentration of Mesh Points near Observation Locations} \label{sec:concen}
Using metric tensors to define a common mesh gives some amount of control over
the location of the mesh points. Specifically, the common mesh can be concentrated near specific points of interest, such as observation locations. In the context of data assimilation, it may be beneficial to have more mesh points near the observation locations $\{x_j^o(t)\}_{j=1}^{N_o}$, where $N_o$ is the number of observations. The locations of the observations do not need to be fixed, as long as the location is known at the observational time $t$.

One strategy for ensuring that there are common mesh points near the observation locations is to simply fix the observation locations as points in the common mesh. However, this can lead to meshes that are ill-conditioned when the mesh points are not allowed to move freely with the dynamics of the solution. Instead of fixing the observation location as a point in the common mesh, a monitor function $\M_K^O$ can concentrate the mesh
near the observation locations. Ideally, $\M_K^O$ is comparable to $\M_K^m$ in some measure at the locations where the mesh requires more points, and quickly decays away from these locations. One such choice is
\[\M_K^O = \sum_{j=1}^{N_o}\chi\left(\|x-x_j^o(t)\|\right)I,\]
where
\[\chi(w) = \left[e^{4w^2}-1+\frac{1}{\max\limits_{K} \sqrt{\det(\M_K^m)}}\right]^{-1}.\]

\begin{equation}\label{eq:concen}
\M_K = \M_K^m \cap \M_K^O .
\end{equation}
If the ensemble meshes alone determine the common mesh, set $\M_K = \M_K^m$. 
Conversely, if the locations of the observations alone determine the common mesh, set $\M_K = \M_K^O$. 
Note that this last scenario is similar to \cite{Du}, where
a fixed, non-uniform common mesh is employed based on the observation locations.


\subsection{Adaptive Localization} \label{sec:adaptive_loc}
A localization scheme can ensure that observations only affect nearby mesh points. There are broadly two categories of localization schemes: domain localization and covariance localization. For a covariance localization scheme, the covariance matrix is adjusted so that the analysis update is less affected by observations that are farther away. Several works explore this type of R-localization scheme, including in an adaptive sense \cite{Wang2018,Moosavi2018,Popov2019}.
The second category of localization schemes is domain localization, where the domain is decomposed into several subdomains. Domain localization ensures that an observation only affects the solution at mesh points within the same subdomain. It will not affect the solution outside of the given subdomain.
One common domain localization scheme is that used in \cite{Hunt}.


We define a metric tensor localization scheme (MT localization) as a domain localization scheme of \cite{Hunt}
(cf. Section~\ref{sec:details-localization}) but with the localization radii calculated for all mesh node as a function of the determinant of the monitor function as follows.
Instead of having one predetermined radius of localization $r$, at each timestep we compute the localization radius for each node:
\begin{equation}\label{eq:dynloc}
r_i = L\;  e^{-\frac{d_i}{2 d_{min}}} , \quad i = 1, ..., M
\end{equation}
where $d_i = \min(\det(\M_K^m(x_i)),c)$, $d_{min} = \min_i d_i$, and $c > 0$ and $L > 0$ are the parameters offering some control over the localization regime. It is not difficult to see that
\begin{equation} \label{eq:bounds}
L\; e^{-\frac{c}{2 d_{min}}}\leq r_i \leq \frac{L}{\sqrt{e}}.
\end{equation}
This shows that the larger the value of $L$, the larger the localization radius can be. A smaller cutoff value $c$ will increase the lower bound of the localization radius, ensuring that localization can still happen, even given a sharp front.

When combining the concentration of mesh points with the MT localization scheme, the localization radius should be calculated solely based on the meshes from the previous time-step and the ensemble solutions. This is done before the common mesh is concentrated near the observations. If not, the localization radius would be incorrectly computed from the concentration scheme instead of the ensemble solutions to the PDE.

\subsection{Remeshing ensemble members}
At the beginning of each observational timestep, the ensemble members $u^{e_i}$ reside on their corresponding meshes $\T_h^{e_i}$. The analysis is computed on the common mesh, and the updated ensemble members are interpolated back to their individual meshes. 
Sometimes the analysis update is enough of a perturbation from the forecast that the meshes that worked for the forecast are no longer suitable for the analysis meshes. 

If the perturbation is large enough, the forecast meshes might not still be appropriate for the analysis. One option is to add extra smoothing cycles before integrating. If necessary, we can remesh to find a mesh suitable for the analysis ensemble. The new meshes for ensemble members must satisfy the equidistribution and alignment criteria for the updated ensemble, and the result is the analysis ensemble $u^{e_i,a}$ residing on the updated meshes $\T_h^{e_i,a}$.


\subsection{Algorithm for DA on Adaptive Meshes}
The technique that has been developed here for DA with an adaptive moving mesh is summarized in Algorithm \ref{DAMMAlg}.
Here we assume that all meshes have the same number of elements and nodes and the same connectivity. A main advantage of this is that those meshes can be viewed as deformations to each other and conservative interpolation schemes can be developed relatively easier between those meshes; e.g., see \cite{Huang2021} or Section~\ref{sec:details-dginterp}. It is important in the DA computation to conserve ensemble members at the interpolation steps 7 and 8 in Algorithm \ref{DAMMAlg} since, otherwise, the mean zero assumption for the model error in (\ref{eq:dtime_mod}) will be violated.

We note that it is not a requirement for meshes to have the same number of elements and nodes and the same connectivity for the metric tensor approach to ensemble DA with adaptive meshing. This approach can also be used with a meshing scheme where mesh points may be added or eliminated based on some pre-defined criteria and/or ensemble meshes and the common mesh can have different numbers of elements and nodes. However, precaution should be taken for interpolation between ensemble meshes and the common mesh to ensure conservation of the ensemble members and therefore mean zero of the model error.

\begin{algorithm}
\caption{EnKF on Adaptive Moving Meshes}\label{alg:MMPDEDA}
\begin{algorithmic}[1]
\Procedure{DA Update on Moving Mesh}{}
\State Compute $\M^{m}$
\State Compute $\M^{O}$
\State Take $\M = \M^{m}, \M^{O}$ or $\M^{m}\cap \M^{O}$
\State Compute common mesh $X$
\State (Optional) Adapt localization scheme based on $\M^{m}$
\State Interpolate $u^{m}$ and observations to common mesh $X$
\State DA update on common mesh $X$
\State Interpolate $u^{m,analysis}$ to individual meshes
\State (Optional) Remesh ensemble meshes
\State Integrate solutions forward in time until next observation
\EndProcedure
\end{algorithmic}\label{DAMMAlg}
\end{algorithm}

\section{Implementation Details}\label{sec:details}
\subsection{Defining Mesh Movement}
\label{sec:details-mm}
The following is a summary of the moving mesh PDE (MMPDE) approach developed in \cite{Huang2006,Huang2015,Huang1994-1,Huang1994-2,Huang2010}. The central idea of the MMPDE moving mesh method is to view any nonuniform mesh as a uniform one in some metric $\M$; that is, the elements have a constant volume and are equilateral in the metric $\M$. These conditions are called the equidistribution and alignment criteria; see (\ref{eq:equidistribution}) and (\ref{eq:alignment}). It has been shown that if a mesh begins as nonsingular (that is, the elements have positive volume), it will remain nonsingular for all time under the MMPDE method. Furthermore, the height of each of the elements will be bounded above and below by some positive constants that depend on $\M$ and the initial mesh \cite[Theorem 4.1]{Huang2018}.

The metric tensor $\M = \M({x})$ is used to control the size, shape, and orientation of the elements of the mesh to be generated. Various metric tensors have been developed in \cite{Huang2003}. For this paper, we consider a Hessian-based metric tensor defined for each element $K \in \T_h$ as
\begin{equation}\label{mer}
   \M_{K} =\det \left(I+\frac{1}{\alpha_h}|H_K(u)|
        \right)^{-\frac{1}{d+4}}
   \left(I+\frac{1}{\alpha_h}|H_K(u)|\right) ,
\end{equation}
where $H_K(u)$ is the Hessian or a recovered Hessian of the state vector $u \in \R^d$ on the element $K$; $|H_K(u)| = Q
\text{diag}(|\lambda_1|,...,|\lambda_d|)Q^T$ with $Q \text{diag}(\lambda_1,
...,\lambda_d)Q^T$ being the eigen-decomposition of $H_K(u)$; and $\alpha_h$
 is a regularization parameter defined through the following algebraic equation:
\[
\sum_{K\in\T_h}|K|\, \det\left(I+\frac{1}{\alpha_h}|H_K(u)|\right)^{\frac{2}{d+4}}
    =2\sum_{K\in\T_h}|K|\,
    \det\left(|H_K(u)|\right)^{\frac{2}{d+4}}.
\]
The metric tensor given by equation \eqref{mer} is known to be optimal for the $L^2$-norm of linear interpolation error \cite{Huang2003}.

To generate the $\M$-uniform mesh $\T_h$, we use here an approach different from (\ref{eq:mesheq}) where the coordinators of the mesh nodes are evolved directly. Instead, we evolve the coordinators of the nodes of an intermediate mesh and obtain the new mesh through this intermediate mesh and interpolation. An advantage of this approach is that its formulation is simpler than that of the direct approach (cf. \cite{Huang2015}). To start with, we introduce the reference computational mesh $\hat{\T}_c$ which is uniform in the Euclidean metric, and the computational mesh $\T_c$. The reference computational mesh $\hat{\T}_c =\{{\hat{\xi}}_j\}_{j=1}^{M}$ has the same connectivity and the same number of vertices and elements as $\T_h$ and stays fixed in the computation. The computational mesh $\T_c =\{{\xi}_j\}_{j=1}^{M}$ serves as an intermediate variable. In this setting, for any element $K$ in $\T_h$, there is a unique corresponding element $K_c$ in $\T_c$. The affine map between $K_c$ and $K$ and its Jacobian matrix are denoted by $F_K$ and $F_K'$, respectively. With this new setting, the equidistribution and alignment criteria for $\M$-uniform meshes have a similar form as those in equations \eqref{eq:equidistribution} and \eqref{eq:alignment}.
\cmt{
It is known that an $\M$-uniform
mesh $\mathcal{T}_h$ satisfies the so-called equidistribution and alignment conditions:
\begin{align}
\label{eq:econd}
& |K|\sqrt{\det(\widetilde{\mathbb{M}}_K)}=\frac{\sigma_h|K_c|}{|\Omega_c|},
&\quad \forall K\in \mathcal{T}_h
\\
\label{eq:alcond}
& \frac{1}{d}\hbox{tr}\left( (F'_K)^{-1}\widetilde{\mathbb{M}}_K^{-1}(F'_K)^{-T}\right)
= \hbox{det}\left(
(F'_K)^{-1}\widetilde{\mathbb{M}}_K^{-1}(F'_K)^{-T}\right)^{\frac{1}{d}},
&\quad \forall K \in \mathcal{T}_h
\end{align}
where $F'_K$ is the Jacobian matrix of the affine mapping:
$F_K: K_c\in \mathcal{T}_c \rightarrow K \in \mathcal{T}_h$,
$\widetilde{\mathbb{M}}_K$ is the average of $\widetilde{\mathbb{M}}$ over $K$,
$\hbox{tr}(\cdot)$ denotes the trace of a matrix, and
\[
|\Omega_c|=\sum\limits_{K_c\in\mathcal{T}_c}|K_c|,\quad \sigma_h
=\sum\limits_{K\in\mathcal{T}_h}|K|\hbox{det}(\widetilde{\mathbb{M}}_K)^{\frac{1}{2}} .
\]
} 
To generate a mesh satisfying these conditions as closely as possible, define an energy functional as
\begin{equation}\label{energy}
\begin{split}
I_h(\mathcal{T}_h,\mathcal{T}_c)
=&\frac{1}{3}\sum_{K\in\mathcal{T}_h}|K|\hbox{det}({\mathbb{M}}_K)^{\frac{1}{2}}\left(
\hbox{tr}((F'_K)^{-1}{\mathbb{M}}^{-1}_K(F'_K)^{-T})\right)^{\frac{3 d}{4}}
\\&
+\frac{1}{3} d^{\frac{3 d}{4}}\sum_{K\in\mathcal{T}_h}|K|\hbox{det}
({\mathbb{M}}_K)^{\frac{1}{2}}
\left (\hbox{det}(F'_K) \hbox{det}({\mathbb{M}}_K)^{\frac{1}{2}} \right )^{-\frac{3}{2}},
\end{split}
\end{equation}
which is a Riemann sum of a continuous functional developed in \cite{Huang2001}
based on mesh equidistribution and alignment.
Note that $I_h(\mathcal{T}_h,\mathcal{T}_c)$ is a function of the coordinates of the nodes of $\mathcal{T}_h$
and $\mathcal{T}_c$.

Taking $\mathcal{T}_h$ as the current mesh $\mathcal{T}^{m}_h$, the MMPDE approach defines the mesh equation as a gradient system of
$I_h(\mathcal{T}_h,\mathcal{T}_c)$,
\begin{equation}\label{MM}
\begin{split}
\frac{d {\xi}_j }{dt}
=-\frac{\hbox{det}({\mathbb{M}}({x_j}))^{\frac{1}{2}} }{\tau}
\left(\frac{\partial I_h(\mathcal{T}^{m}_h,\mathcal{T}_c)}{\partial {\xi}_j}\right)^T,
\quad j=1,\cdots,M
\end{split}
\end{equation}
where ${\partial I_h }/{\partial {\xi}_j}$ is considered as a row vector,
$\tau>0$ is a parameter used to adjust the response time of mesh movement to the
changes in ${\mathbb{M}}$.

Define the function $G$ associated with the energy \eqref{energy} as
\begin{equation}\label{G}
G(\mathbb{J},\hbox{det}(\mathbb{J}))=
\frac{1}{3}\hbox{det}({\mathbb{M}}_{K})^{\frac{1}{2}}
(\hbox{tr}(\mathbb{J}{\mathbb{M}}_{K}^{-1}\mathbb{J}^T))^{\frac{3 d}{4}}
+\frac{1}{3} d^{\frac{3 d}{4}}\hbox{det}({\mathbb{M}}_{K})^{\frac{1}{2}}
\left (\frac{\hbox{det}(\mathbb{J})}{\hbox{det}({\mathbb{M}}_{K})^{\frac{1}{2}}}
\right )^{\frac{3}{2}},
\end{equation}
where $\mathbb{J}=(F'_K)^{-1} = E_{K_c}E_K^{-1}$ and the edge matrices of $K$ and $K_c$ are
$E_K=[{x}_1^K-{x}_0^K,\;\cdots,\; {x}_d^K-{x}_0^K]$ and
$E_{K_c}=[{\xi}_1^K - {\xi}_0^K,\;\cdots,\; {\xi}_d^K -{\xi}_0^K ]$, respectively.
Using the notion of scalar-by-matrix differentiation \cite{Huang2015}, it is not difficult
to find the derivatives of $G$ with respect to $\mathbb{J}$ and $\hbox{det}(\mathbb{J})$ as
\begin{align}
&\frac{\partial G}{\partial\mathbb{J}}=
\frac{d}{2}\det(\mathbb{M}_{K})^{\frac{1}{2}}
(\hbox{tr}(\mathbb{J}{\mathbb{M}}_K^{-1}\mathbb{J}^T))^{\frac{3 d}{4}-1}
{\mathbb{M}}_K^{-1}\mathbb{J}^T,
\label{partial-J}
\\
&\frac{\partial G}{\partial \hbox{det}(\mathbb{J})}=
\frac{1}{2} d^\frac{3 d}{4}\det({\mathbb{M}}_{K})^{-\frac{1}{4}}\det(\mathbb{J})^{\frac{1}{2}}
\label{partial-detJ}.
\end{align}
Substituting \eqref{G}-\eqref{partial-detJ} into \eqref{MM} yields
\begin{equation}\label{xim}
\begin{split}
\frac{d{\xi}_j}{dt}=
\frac{\hbox{det}({\mathbb{M}}({x_j}))^{\frac{1}{2}}}{\tau}
\sum_{K\in\omega_j}|K|{v}^K_{j_K},
\quad j=1,\cdots,M
\end{split}
\end{equation}
where $\omega_j$ is the element patch associated with the vertex ${x}_j$, $j_K$ is the
local index of ${x}_j$ on $K$, and ${v}^K_{j_K}$ is the local velocity contributed
by the element $K$ to the vertex $j_K$. The local velocities ${v}^K_{j_K},\; j_K=1,
\cdots,d$ are given by
\begin{equation}\label{vjk}
\begin{split}
\left[
  \begin{array}{c}
    ( {v}_1^K  )^T  \\
    ( {v}_2^K  )^T   \\
    \vdots\\
    ( {v}_d^K  )^T   \\
   \end{array}
 \right]
=
- E_K^{-1}\frac{\partial G }{\partial \mathbb{J} }
- \frac{\partial G}{\partial \det(\mathbb{J})}\frac{ \det( E_{K_c} )}
{\det(E_K)} E_{K_c}^{-1},
\quad
{v}^K_0=-\sum_{j_K=1}^d{v}^K_{j_K} .
\end{split}
\end{equation}

Integrating the mesh equation \eqref{xim} over a physical time step, with the proper modifications for boundary vertices and with the initial mesh $\hat{\mathcal{T}}_c$, yields the new computational mesh $\T^{new}_c$ which forms a correspondence with the current mesh $\T_h^m$, i.e., $\T_h^m = \Phi_h(\mathcal{T}_c)$. The new common mesh $\mathcal{T}_h^{new}$ is defined as $\mathcal{T}_h^{new}=\Phi_h(\hat{\mathcal{T}}_c)$, which can be computed using linear interpolation.

It is common practice in moving mesh computation to smooth the metric tensor for smoother meshes. To this end, we apply a low-pass filter \cite{Huang2010} to the metric tensor several sweeps every time it is computed.

\subsection{DG Discretization}
\label{sec:details-dg}

Numerical results will be presented in Section~\ref{sec:results} to illustrate the DA procedure using the 1D and 2D inviscid Burgers equations. Since the Burgers equation is hyperbolic and can have discontinuous solutions such as shocks, we use the discontinuous Galerkin method (DG) for its spatial discretization. DG is a type of finite element method with the trial and test function spaces consisting of discontinuous, piecewise polynomials. It is known to be a particularly powerful numerical tool for the simulation of hyperbolic problems and has the advantages of high-order accuracy, local conservation, geometric flexibility, suitability for handling mesh-adaptivity, extremely local data structure, high parallel efficiency, and a good theoretical foundation for stability and error estimates. Over the last few decades, the DG method has been used widely in scientific and engineering computation.

Specifically, we use a quasi-Lagrangian moving mesh DG method (MMDG) to solve the Burgers equation on moving meshes \cite{Huang2021}. The method treats the mesh movement continuous in time and leads to an extra convective term (cf. equation (\ref{eq:mmpde})) in the resulting discrete equations to reflect the mesh movement. Importantly, the method is (mass) conservative so that the model error has mean 0. In our computation, we use piecewise linear polynomials (P$^1$-DG) for spatial discretization and a third-order Strong Stability Preserving (SSP) Runge-Kutta scheme for temporal discretization. The reader is referred to \cite{Huang2021} for details of the MMDG method.

\subsection{DG Interpolation}
\label{sec:details-dginterp}

From Algorithm~\ref{alg:MMPDEDA}, we see that interpolation is needed between the ensemble meshes and the common mesh. Since these meshes are assumed to be deformations to each other, we can perform interpolation by solving a differential equation.

To be specific, we consider the deforming meshes $\T_h^{old}$ and $\T_h^{new}$ and the state variable $u$ that needs to be interpolated. We define a deforming mesh $\T_h(\zeta)$ (with $\zeta \in [0,1]$) from $\T_h^{old}$ to $\T_h^{new}$ as a mesh with the nodal positions and velocities as
\begin{align}
    x_i(\zeta) &= (1-\zeta)x_i^{old} + \zeta x_i^{new}, \quad i = 1,\dots, M\\
    \dot x_i &= x_i^{new} - x_i^{old}, \quad i = 1, \dots, M.
\end{align}
Then the interpolation can be viewed as solving the following linear convective PDE on the moving mesh $\T_h(\zeta)$ over $\zeta \in [0,1]$:
\begin{equation}\label{eq:DGeqn}
    \frac{\partial u}{\partial \zeta}(x,\zeta) = 0, \quad (x,\zeta)\in \Omega\times(0,1],
\end{equation}
with the initial values $u(x,0)$ as those of $u$ on $\T_h^{old}$.

A DG-interpolation scheme has been studied in \cite{Huang2021} where DG and a SSP Runge-Kutta scheme are used to discretized (\ref{eq:DGeqn}) in space and time, respectively. The scheme is conservative and positivity-preserving and can be high-order. As mentioned before, the mass conservation is important for interpolation between ensemble meshes and the common mesh to maintain the mean zero feature of the model error.
Like for the PDE solver, we use this scheme for interpolation with piecewise linear polynomials (P$^1$-DG)
for spatial discretization and a third-order SSP Runge-Kutta scheme for temporal discretization; see \cite{Huang2021} for detail.

\subsection{DA Implementation}
For our numerical experiments we employ a Local Ensemble Transform Kalman Filter (LETKF) based on \cite{Hunt}. In general, an ETKF code uses a linear transform to have control over the resulting sample covariance so that the covariance of the analysis update, $P^f$, exactly satisfies the identity $P^f_{j+1} = (I-K_{j+1}H)P^b_{j+1}$.

Let $\hat m_{n+1}$ be the ensemble mean at time $t_{n+1}$, and let $\hat u_{n+1}^{(i)}$, $i = 1, \dots N_e$ be the ensemble forecast. Then define the perturbation matrix
\begin{equation}
    \hat X_{n+1} = \frac{1}{\sqrt{N_e -1}} \left[\hat u_{n+1}^{(1)} - \hat m_{n+1}, \hat u_{n+1}^{(2)} - \hat m_{n+1} , \dots, \hat u_{n+1}^{(N_e)} - \hat m_{n+1}\right],
\end{equation}
and the sample covariance is given by $P^b_{n+1} = \hat X_{n+1} \hat X_{n+1}^T$.
Consider the following transformation:
\begin{equation} T_{n+1} = \left[ I + \left(H\hat X_{n+1}\right)^T R^{-1}\left(H\hat X_{n+1}\right)\right]^{-1},
\end{equation}
and set $X_{n+1} = \hat X_{n+1}T_{n+1}^\frac{1}{2}$. Then
\begin{equation}
    P^f_{n+1} \coloneqq X_{n+1}X_{n+1}^T = (I-K_{n+1}H)P^b_{n+1},
\end{equation}
as desired.




\subsection{Localization Techniques}
\label{sec:details-localization}

We will compare the metric tensor based localization scheme developed in Section \ref{sec:adaptive_loc} with some commonly used localization schemes.
Recall from \eqref{eq:DA_analysis} the Kalman gain $\mathbf{K}_{n+1}$.
R-localization modifies the Kalman gain through the Schur product of a localization matrix $\rho$ with the covariance matrix. One of the most common ways to define $\rho$ is through the Gaspari-Cohn (GC) correlation function, which is a fifth order polynomial that decays to zero.
\begin{equation} \label{eq:GC_poly}
\mathcal{C}(r) = \begin{cases}
-\frac{1}{4}r^5 + \frac{1}{2}r^4 + \frac{5}{8}r^3 - \frac{5}{3}r^2+1, & 0 \leq r \leq 1\\
\frac{1}{12}r^5 - \frac{1}{2}r^4 + \frac{5}{8}r^3 + \frac{5}{3}r^2 - 5r + 4 - \frac{2}{3}r^{-1}, & 1 < r \leq 2\\
0, & 2<r .
\end{cases}
\end{equation}
GC localization can be applied to either the model space or the observation space. In the model space, the localization matrix is given by
\begin{equation} \label{eq:GC_mod}
    \rho_{ij} = \mathcal{C}(|x_i-x_j|/L),\quad i, j = 1, ..., M
\end{equation}
where $x_i$ and $x_j$ are the positions of the $i$th and $j$th nodes and $L > 0$ is a pre-determined localization radius. This matrix is Schur-multiplied with the ensemble covariance matrix, resulting in the Kalman update
\begin{equation}
    K_{GCmod} = \left(\rho \circ P^b\right)H^T\left(H\left(\rho\circ P^b\right)H^T + R\right)^{-1}.
\end{equation}

For localization in the observation space, the localization matrix must be Schur multiplied to both $\left(P^b H^T\right)$ and $\left(HP^b H^T\right)$. Therefore, two localization matrices are needed,
\begin{align}
    \left[\rho_1\right]_{i,j} &= \mathcal{C}(|x_i-y_j|/L), \quad \quad i, j = 1, ..., M
     \label{eq:GC_obs1}\\
    \left[\rho_2\right]_{i,j} &= \mathcal{C}(|y_i-y_j|/L), \quad \quad i, j = 1, ..., M.
    \label{eq:GC_obs2}
\end{align}
Then the Kalman gain is given by
\begin{equation}
    K_{GCobs} = \rho_1 \circ \left(P^b H^T\right)\left(\rho_2\circ \left( HP^b H^T \right)+ R\right)^{-1}.
\end{equation}

We next outline the domain localization scheme developed in \cite{Hunt}.
 Define $r$ to be the radius of localization. Then for every mesh point $x_i$, if there is an observation $y$ located at the point $x^o$ such that $\|x_i - x^o\|_2 \leq r$, the analysis is given by the EnKF update in equation \eqref{eq:DA_analysis}. If not, then the analysis update is equal to the forecast predicted by the model.

In many implementations, this radius of influence is predetermined and constant over time and space. However, there may be instances where this should be dynamic in time and space. For example, if the solution has a traveling front or shock that travels across the domain, there should be a smaller radius of localization near the region where the gradient is large and a larger radius where the solution is relatively constant, and this localization scheme should move across the domain with the traveling front or shock.

To see this, consider a solution $u_n(x) \in \R$ with a single observation at $x^o$ and a large gradient beginning at $x_k$ just past this observation location; that is, $0 < x_k-x^o \ll 1$. For $x_i$ outside of the localization radius, the observation does not affect the analysis, so $u_{n+1}(x_i) = \hat u_{n+1}(x_i)$. For $x_i$ close to the observation, consider the EnKF update for a single ensemble member, omitting the ensemble superscripts and time subscripts:
\begin{equation}
    u(x_i) = (I - \mathbf{K}H)\hat u(x_i) + \mathbf{K}y,
\end{equation}
where $\mathbf{K} = P^bH^T\left(HP^bH^T + R\right)^{-1}$.
For the sake of simplicity, assume $R = \alpha I$ and $H = e_k^T$,
where $e_k$ denotes the $k^{th}$ unit vector. Then the analysis update at the point $x_i$ close to the observation gives us
\[
u^a(x_i) = u^f(x_i) + \mu P^b_{ik}\quad \text{with}\quad \mu = \left[P^b_{kk} + \frac{1}{\alpha}\right](y-x^f_k) .
\]
In particular, the larger the value of $P^b_{ik}$ is, the greater the difference will be between the analysis and the forecast. This can be especially problematic if $P^b_{ik}$ is large for $i<k$, 
that is, after the shock. In that case, the analysis will be changed to be much closer to the observations, and the steep front will be smoothed out.

This problem is avoided by using a smaller radius of localization near the shock and a larger radius of localization farther away from it. Fortunately, the monitor function obtained in the adaptive moving mesh algorithm determines where mesh points will be closer together and where they will be spread far apart; by proxy, this shows where the shock or other feature of interest exists. In this way, the monitor function can inform what the localization parameter should be over time and space, enabling a dynamic update of the localization variable.

LETKF employs localization within the ensemble transform Kalman filter so that only model variables located at mesh points within a predetermined radius of an observation will assimilate that observation. This not only localizes the influence of observations but also provides a dimension reduction by creating reduced dimensional subproblems on which assimilation is performed independently.

A covariance-based localization scheme, which uses a Schur product applied to the covariance matrix, is problematic when working in this reduced dimension. For example, the reduced dimension implementation uses $\hat X_{n+1}$, taken from the Cholesky decomposition of the covariance matrix $P^b_{n+1}$, rather than $P^b_{n+1}$ itself. Covariance localization would use a Schur product to adjust this covariance matrix, but in doing so, it could result in a matrix with negative eigenvalues. Therefore, a covariance localization scheme is not easily implemented into the LETKF code. For the experiments where we compare the MT localization scheme to covariance localization schemes like Gaspari-Cohn, we use a traditional ETKF code.

\section{Numerical Results}\label{sec:results}

The following presents the application of these methods to the one and two dimensional inviscid Burgers equations. We generate synthetic observations by sampling from a truth run, obtained by solving this equation on an adaptive moving mesh. The ensemble members are initialized as perturbations of the initial conditions. Efficacy of the DA scheme is measured by the root mean squared error (RMSE), which is calculated as
\begin{equation}
    \text{RMSE} \coloneqq \frac{1}{\sqrt{M}} \|u^{\text{truth}} - \bar u\|_2,
\end{equation}
where $\bar u$ is the analysis mean. A DA procedure is generally considered stable if its asymptotic behavior is on the order of the square root of the norm of the observation error. The RMSE in the experiments that follow is averaged over 10 runs.

\subsection{Common Experimental Set-Up}
The next two sections explore the use of the LETKF on the one- and two-dimensional inviscid Burgers equation.
We perform identical twin experiments where the truth is generated using no model noise and observations are formed from the truth by applying the observation operator and adding noise using the observation error model with covariance matrix $R = 0.01I$.
Among the parameters to be tuned are the number of mesh points and the number of ensemble members. Generally speaking, more mesh points correspond to more accurate numerical solutions, lowering the RMSE. However, moving mesh methods generally require fewer mesh points than a fixed, uniform mesh. For the 1D Burgers experiments, we use $M=50$ mesh points, and for the 2D Burgers experiments, we use a $M = 15\times 15$ mesh. Increasing the number of mesh points beyond the values chosen had little impact on the RMSE. The rule of thumb for ensemble-based DA schemes is that the number of ensemble members should be large enough to span the unstable subspace. For both the 1D and 2D experiments, we found $N_e = 5$ ensemble members to be sufficient. That is, for larger $N_e$, we found that there was no substantial improvement in RMSE.


 In the 1D Burgers experiments we observe the truth at $D = 5$ locations with an observational timestep of $\Delta t = 0.5$, and then add artificial observation error ($\eta \sim \mathcal{N}(0,0.01)$). To avoid the observations all occurring in one region of the spatial domain, we space them linearly throughout the domain, and then perturb the locations by a small amount. These observation locations are randomly chosen, but once chosen at the beginning of each trial, they remain fixed. The truth (and observations) are taken from a fine mesh (100 mesh points) to ensure that it is fully resolved. For the 2D Burgers experiments, we use an observational timestep of $\Delta t = 0.5$ with $D = 16$ observation locations that are also uniformly spaced in a $4\times 4$ mesh, and then perturbed. Unless otherwise stated, the observational error covariance $R$, model error covariance $\Sigma$, and prior distribution $P_b^0$ are all set to $0.01I$.

For both the 1D and the 2D cases, the localization radius and inflation factors are tuned simultaneously. For MT localization, this involves determining the parameter $L$ shown in equation \eqref{eq:dynloc}, which directly controls the maximum radius of localization. In the 1D case, we choose not to artificially limit the minimum value of the MT localization radius by choosing $c$ larger than the maximum of $\M_K$ in \eqref{eq:dynloc}; numerical results suggest that $c=8$ is sufficient. In the 2D case, we keep $c=8$ and note that this does affect the minimum localization radius, but that it also results in a stable DA scheme.
For the GC localization schemes, this tuning experiment involves tuning the parameter $L$ as given in equations \eqref{eq:GC_mod}, \eqref{eq:GC_obs1}, and \eqref{eq:GC_obs2}. After each of the localization schemes has been tuned, the time series RMSE of the different localization schemes is directly compared. Finally, we compare results when choosing $\M^{m}$, $\M^{O}$, or $\M^{m}\cap\M^{O}$ for the common mesh. When considering long-run RMSE results, we consider the RMSE only after the DA scheme has stabilized so as to evaluate the asymptotic behavior of the DA procedure. Based on numeric results, we present the RMSE for 1D Burgers on the time interval $[25,100]$ and for 2D Burgers on the time interval $[15, 50]$. These experiments and the parameters used can be found in Table \ref{tab:parameters1D}.

   \begin{center}
        \begin{table}
       \begin{tabular}{|c|p{26mm}|c|c|p{14mm}|c|p{15mm}|c|}
    \hline
        \textbf{Experiment} & \textbf{Description} & \textbf{Model} & $\mathbf{\Sigma}, \mathbf{R}$ &  \textbf{Inflation Factor} & $\mathbf{M}$ &  \textbf{Loc. Scheme} & \textbf{Mesh Choice}\\
        \hline
        1 & Localization and Inflation Tuning & \eqref{eq:Burgers}& $0.01\cdot I$ & varies & 50 &  varies &  $\M^{m}\cap\M^O$\\
        \hline
        2 & Compare Loc. schemes &  \eqref{eq:Burgers} & $0.01\cdot I$ & varies & 50 & varies & $\M^{m}\cap\M^O$\\
        \hline
        3 & Compare Meshes & \eqref{eq:Burgers}& $0.01\cdot I$ & 1.1 & 50 & MT  & varies\\
        \hline
        4 & Compare Errors & \eqref{eq:Burgers}& varies & 1.1 & 50 & MT & $\M^{m}\cap \M^O$\\
        \hline
        \hline
         5 & Localization and Inflation Tuning & \eqref{eq:Burgers2d} & $0.01\cdot I$ & varies & 225 &  varies &  $\M^{m}\cap\M^O$\\
        \hline
        6 & Compare Loc. schemes & \eqref{eq:Burgers2d} & $0.01\cdot I$ & varies & 225 & varies & $\M^{m}\cap\M^O$\\
        \hline
        7 & Compare Meshes & \eqref{eq:Burgers2d} & $0.01\cdot I$ & 1.1 & 225 & MT  & varies\\
        \hline
        8 & Compare Errors & \eqref{eq:Burgers2d} & varies & 1.1 & 225 & MT & $\M^{m}\cap \M^O$\\
        \hline
        9 & Noisy Data & \eqref{eq:Burgers2d} & $0.01\cdot I$; data varies & 1.1 & 225 & MT & $\M^{m}\cap \M^O$\\
        \hline
        10 & Interpolation & \eqref{eq:Burgers2d} & $0.01\cdot I$ & 1.1 & 225 & MT & $\M^m\cap\M^O$\\
        \hline
    \end{tabular}
    \caption{Summary of experiments.} \label{tab:parameters1D}
    \end{table}
    \end{center}

\subsection{Inviscid Burgers Equation - A One-Dimensional Example}\label{sec:Burgers' Section}
Consider the one-dimensional inviscid Burgers equation as given by
\begin{equation}
\label{eq:Burgers}
\frac{\partial u}{\partial t} + u\frac{\partial u}{\partial x} = 0,
\hspace{5 mm} x \in [0,S), \hspace{5 mm} t \in (0,T],
\end{equation}
with initial condition
\begin{equation}
\label{eq:bgm_ic}
	u(x,0) = \frac{1}{2}+\sin{\left(\frac{2\pi}{S}x\right)}
\end{equation}
and periodic boundary conditions. As time progresses, a shock forms and propagates to the right. This travelling shock makes the 1D inviscid Burgers equation a good candidate for an adaptive moving mesh, as more mesh points are needed near the location of the shock to sufficiently resolve the numerical solution. For the experiments that follow, we consider the spatial domain $[0,20)$ so that $S=20$ over the time interval $(0,100]$.

\subsubsection{Experiment 1: Tuning localization and inflation}
The localization parameters for the various localization schemes are tuned simultaneously with the inflation parameter. Equation \eqref{eq:GC_poly} combined with either \eqref{eq:GC_mod} or \eqref{eq:GC_obs1}-\eqref{eq:GC_obs2} implies that choosing a GC localization parameter of $L=10$ for the GC localization would imply that the entire domain is affected by the localization scheme. Therefore, in this experiment, the GC localization parameters vary from $L = 0.5$ to $L=10$. Similarly, a localization parameter of $L=20$ means that the entire domain could potentially be affected by the MT localization scheme. For the MT localization tuning, $L$ varies from $L=0.5$ to $L=20$. For all three localization schemes, the multiplicative inflation factor varies from $\rho = 1$ (no inflation) to $\rho = 1.5$.
\begin{figure}[ht]
    \centering
    \includegraphics[width = 0.3\textwidth]{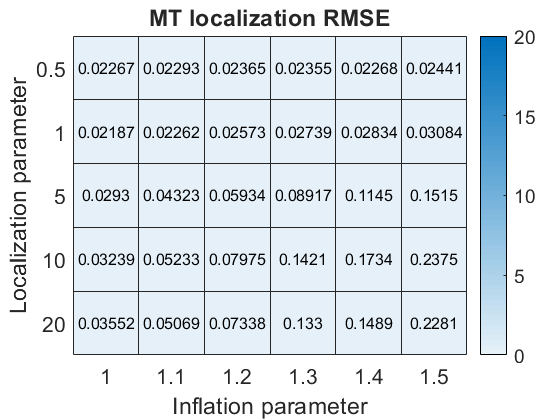}\hfill
    \includegraphics[width=0.3\textwidth]{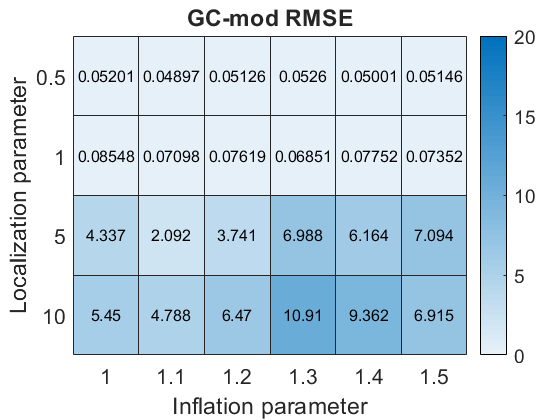}\hfill
    \includegraphics[width=0.3\textwidth]{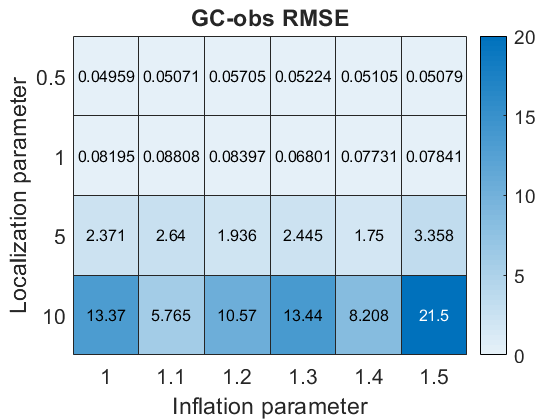}
    \caption{Results for simultaneous tuning of localization and inflation for the MT, GC-mod, and GC-obs localization schemes.}
    \label{fig:tune_rho_irad}
\end{figure}

As seen in Figure \ref{fig:tune_rho_irad}, the MT localization scheme is robust with respect to the tuning of the localization parameter. The GC localization requires much more careful tuning in both the model space and observation space. More specifically, the GC localization schemes work well when the localization parameter $L$ is less than or equal to 1. If the localization parameter is greater than 1, the Schur product of the localization matrix and the covariance matrix is no longer positive definite, and the localization scheme performs poorly.

\begin{figure}[ht]
    \includegraphics[width = .52\textwidth]{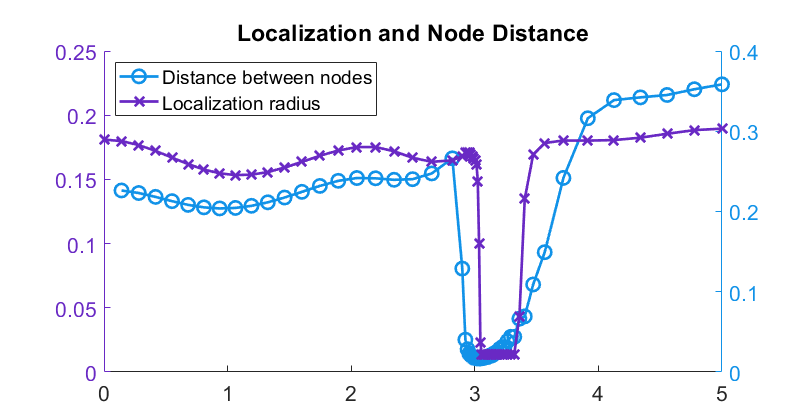}\hfill
    \includegraphics[width=.52\textwidth]{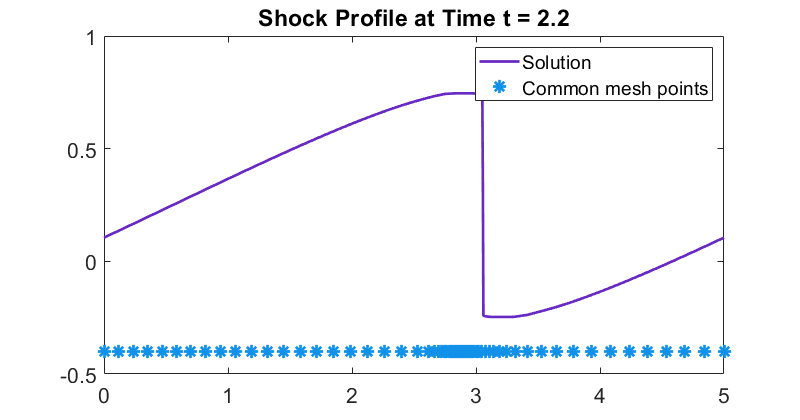}
    \caption{Relationship between nodal distance and adaptive localization parameter. The distance between the nodes is positively correlated with the localization radius.}
    \label{fig:iradrelationship}
\end{figure}

As explained in Section \ref{sec:adaptive_loc}, The MT localization scheme is a domain localization scheme based off $\M_K^m$, the metric tensor intersection of the metric tensors for the ensemble meshes. Since each of the ensemble meshes will have a concentration of points near a shock or large gradient, the adaptive common mesh will also have more mesh points where the solution requires a finer resolution. MT localization reduces the localization radius at these types of interfaces. The relationship between the distance between nodes and the localization variable is shown in Figure \ref{fig:iradrelationship}.

\subsubsection{Experiment 2: Localization Scheme Comparison}
As shown in Experiment 1, the MT localization is much more robust with respect to the tuning parameters than either of the GC localization schemes. Consider the following tuned localization schemes:
\begin{itemize}
\item MT localization with multiplicative inflation parameter $1.1$ and localization parameter $1$.
\item GC-mod localization with multiplicative inflation parameter $1.1$ and localization parameter $0.5$.
\item GC-obs localization with multiplicative inflation parameter $1.0$ and localization parameter $0.5$.
\end{itemize}

\begin{figure}[ht]
    \centering
    \includegraphics[width=.7\textwidth]{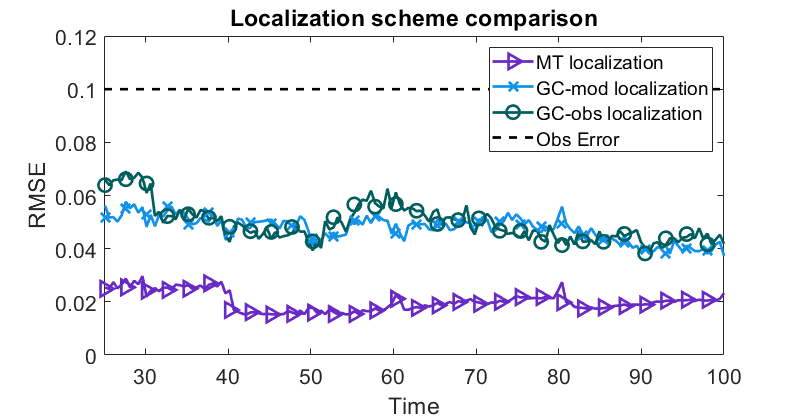}
    \caption{Time series RMSE comparison for the different localization schemes for the 1D inviscid Burgers equation. Each trial uses the tuned localization parameters found in Experiment 1.}
    \label{fig:compare_loc_rmse}
\end{figure}

We compare the time series RMSE of each of these localization schemes in Figure \ref{fig:compare_loc_rmse}. In these experiments, $\Sigma = 0.01I$, so the localization scheme is considered a success if it is on the order of $0.1$. Both the tuned GC-mod localization and the MT localization can be considered a success in this metric, but the MT localization has a considerably lower RMSE than both GC localization schemes.

\subsubsection{Experiment 3: Choice of common mesh}
\begin{figure}[h]
    \centering
    \includegraphics[width = .7\textwidth]{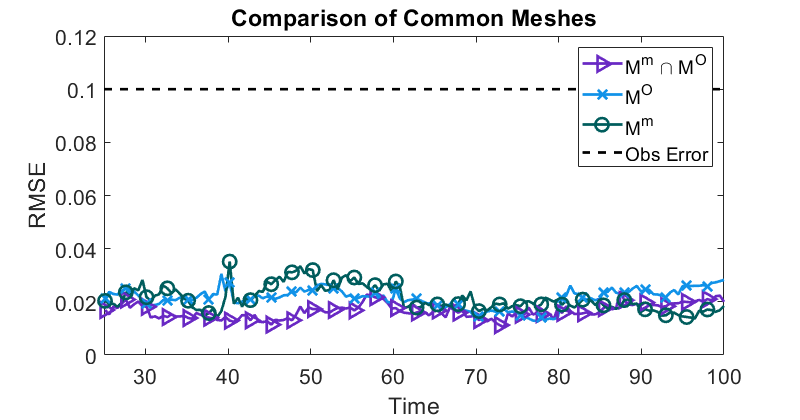}
    \caption{Time series RMSE for tuned DA with the 1D inviscid Burgers equation using different choices for the common mesh.}
    \label{fig:mesh_rmse1D}
\end{figure}

\begin{figure}[ht]

\centering
\includegraphics[width=.3\textwidth]{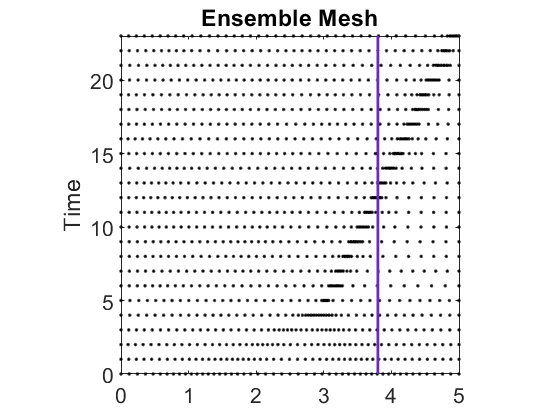}\hfill
\includegraphics[width=.3\textwidth]{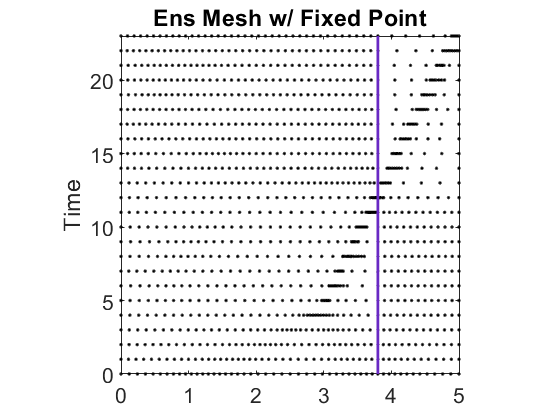}\hfill
\includegraphics[width=.3\textwidth]{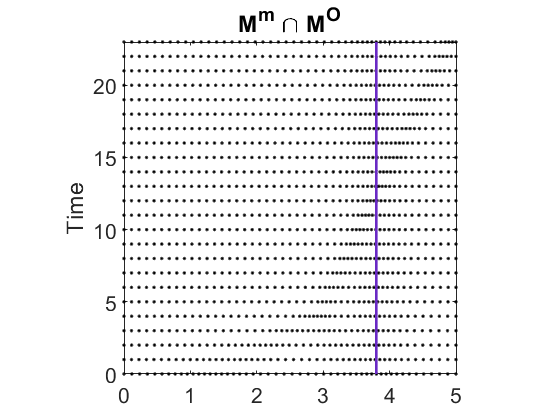}

\caption{A comparison of meshing variations. In the leftmost plot, the mesh evolves naturally with the metric tensor intersection of the ensemble members. In the middle plot, the observation location $x=3.8$ is fixed in the mesh to reduce observation interpolation error. In the right plot, the mesh evolves according to the metric tensor intersection, but it is also concentrated near the location $x=3.8$.}
\label{fig:compareMeshes}

\end{figure}

We consider three choices for the metric tensor associated with the common mesh $\M \in \{\M^{O},\M^{m},\M^{m}\cap\M^{O}\}$. A comparison of these common meshes is shown in Figure \ref{fig:mesh_rmse1D}, using the tuned localization and inflation parameters found in Experiment 1.
While all three choices for the common mesh produce stable DA procedures, using $\M^{O}$ or $\M^{m}\cap \M^{O}$ will concentrate the common mesh near the observation locations, reducing the interpolation error.
 If the interpolation of observations must be avoided at all costs, the user can also specify the observation locations as fixed points in the common mesh. However, this can lead to increased skewness and potential singularity in the mesh, as shown in Figure \ref{fig:compareMeshes}.

The time evolution of the meshes for the 1D inviscid Burgers equation are shown in Figure \ref{fig:compareMeshes}. To better illustrate this phenomenon, we consider a reduced spatial domain of $[0,5]$ $(S=5)$ and update the initial condition accordingly. As the shock forms and propagates to the right, the mesh points of $\M^{m}$ evolve with the shock, as shown in the leftmost plot of Figure \ref{fig:compareMeshes}. However, if an observation location, say at $x = 3.8$, is fixed in the common mesh, it will prohibit the nodes before it from moving past it, as shown in the middle plot of Figure \ref{fig:compareMeshes}. Choosing $\M^{m}\cap\M^{O}$ achieves the goal of following the solution dynamics while still concentrating the mesh near the observation location. (Note that if $\M = \M^{O}$ and the location is static, this common mesh will not change in time.)

\begin{figure}
    \centering
    \includegraphics[width=0.5\textwidth]{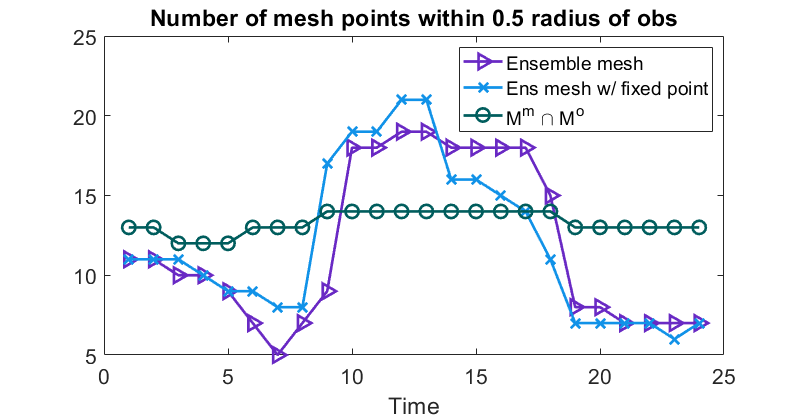}\hfill
    \includegraphics[width=0.5\textwidth]{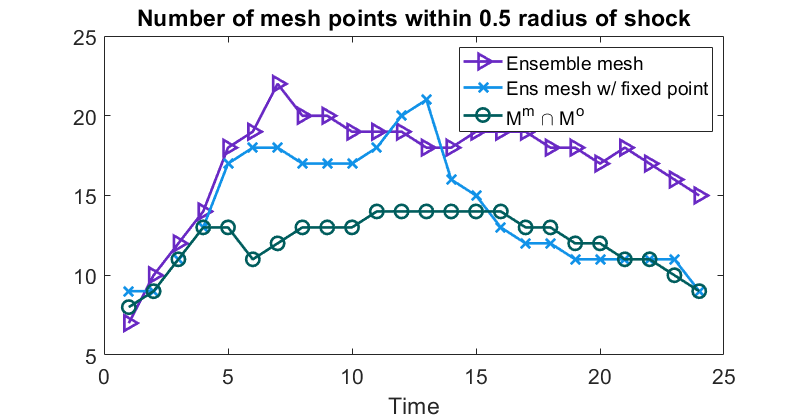}
    \caption{The number of mesh points near the observations varies widely for the ensemble mesh with or without the fixed point. The number of points near the observation has less variance when the mesh points are also concentrated near the observation location.}
    \label{fig:points_count}
\end{figure}

As shown in Figure \ref{fig:compareMeshes}, there are fewer mesh points near the observation in the ensemble mesh with or without the fixed point than there are in the mesh that also concentrates the mesh points near the observation. Figure \ref{fig:points_count} tallies the number of mesh points within the radius $0.5$ of the shock and of the observation for each of the above meshes. The average and variance of these values is given in Table \ref{tab:points_count}. While the meshes formed from the ensemble members, with and without the fixed point, do have an adequate number of mesh points near the observation at times, that only coincides with the shock passing through the observation point. At other times, there are relatively few mesh points near that observation. This is easily seen in the large variances in Table \ref{tab:points_count}.

\begin{table}[ht]
    \centering
    \begin{tabular}{l|c|c|c}
         &  \bf{Ensemble Mesh} & \bf{Ens Mesh w/ Fixed Point} & $\mathbf{\M^{m}\cap \M^O}$\\
         \hline
        \bf{Near Shock} & $(17.13,11.51)$ & $(14.04,13.61)$ & $(12.13,3.16)$\\
        \bf{Near Observation} & $(11.83,25.28)$ & $(11.96,23.87)$ & $(13.29,0.48)$
    \end{tabular}
    \caption{Mean and variance of the number of points near (within a $0.5$ radius of) the shock and observation for each of the three meshes listed above.}
    \label{tab:points_count}
\end{table}

\begin{figure}
    \centering
    \includegraphics[width=0.7\textwidth]{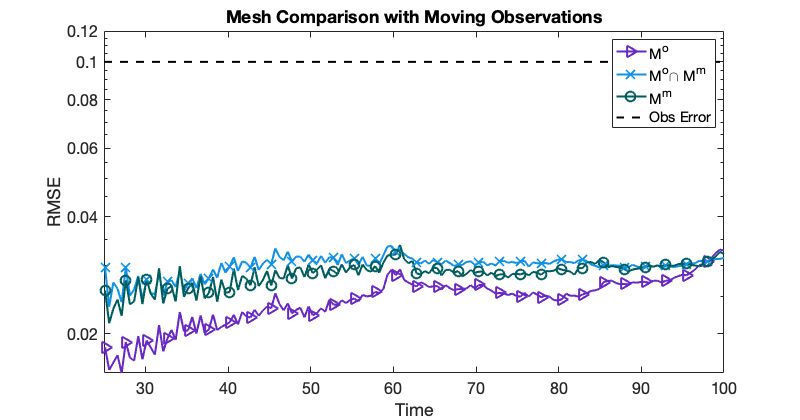}
    \caption{The common mesh based on observation locations performs better than the other choices of common mesh in the case where we have observations that are located away from the shock.}
    \label{fig:moving_obs}
\end{figure}

The benefit of having mesh points concentrated near observations is evident in the case where observations are sparse and occur in places where the mesh points would otherwise not be concentrated. For example, consider the 1D inviscid Burgers equation with two observation locations located away from the shock, and suppose that the location of the observations move at the same rate as the shock as it propagates forward in time. In that case, the observation locations may not see a concentration of mesh points pass through unless the user specifically concentrates the mesh in that area. This is the setup of the experiment shown in Figure \ref{fig:moving_obs}. Over time, the solution flattens due to numerical dissipation, and all choices of the common mesh perform equally well. However, the mesh based on observation locations outperforms the other choices of common mesh before the shock decays.

\cmt{
The full time series of the RMSE for the 1D Burgers using 40 mesh points, inflation value of 1.2, and 5 ensemble members is shown in \ref{fig:1DRMSE}. It uses the common mesh computed by $\M^{O} \cap \M^{m}$, and an adaptive localization scheme with the maximum radius equal to 0.5.
\begin{wrapfigure}{r}{0.3\textwidth}
    \centering
    \includegraphics[width=.3\textwidth]{Fig/1DRMSE.png}
    \caption{Time series RMSE for tuned 1D Burgers.}
    \label{fig:1DRMSE}
\end{wrapfigure}
} 

\subsubsection{Experiment 4: Comparison of error covariances}
\begin{figure}[ht]
    \centering
    \includegraphics[width = .8\textwidth]{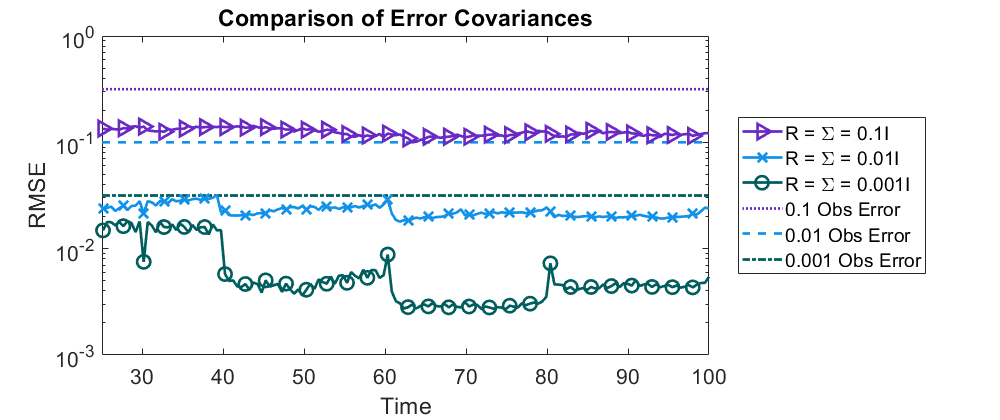}
    \caption{Time series RMSE for tuned MT localization with different error covariances.}
    \label{fig:covariances}
\end{figure}

Using the tuned MT localization scheme from experiment 1, we test the robustness of this DA procedure on \eqref{eq:Burgers} with different error covariances. 
The results in Figure \ref{fig:covariances} show stable RMSE, especially for larger error covariances. The spikes in the RMSE for smaller error covariance correspond to times in which the shock passes through the boundary. In this  implementation the PDE satisfies periodic boundary conditions while the mesh satisfies Dirichlet boundary conditions so that when shock passes through the boundary, the PDE is not approximated as accurately. 


\subsection{Inviscid Burgers Equation - A Two-Dimensional Example}
Consider the two-dimensional Burgers equation
\begin{equation}
    \label{eq:Burgers2d}
    u_t + u u_x + u u_y =0,
\end{equation}
with $\Omega = (-0.5,1)\times (-0.5,1)$, $t \in (0,5]$, and periodic boundary conditions. Given the initial condition
\[
u = \exp(-\gamma(x^2 + y^2)),
\]
with $\gamma = -\log(10^{-16})$, the solution will
have a Gaussian bump that propagates diagonally to the upper right corner of the domain.

\begin{figure}
     \centering
     \includegraphics[width=.24\textwidth]{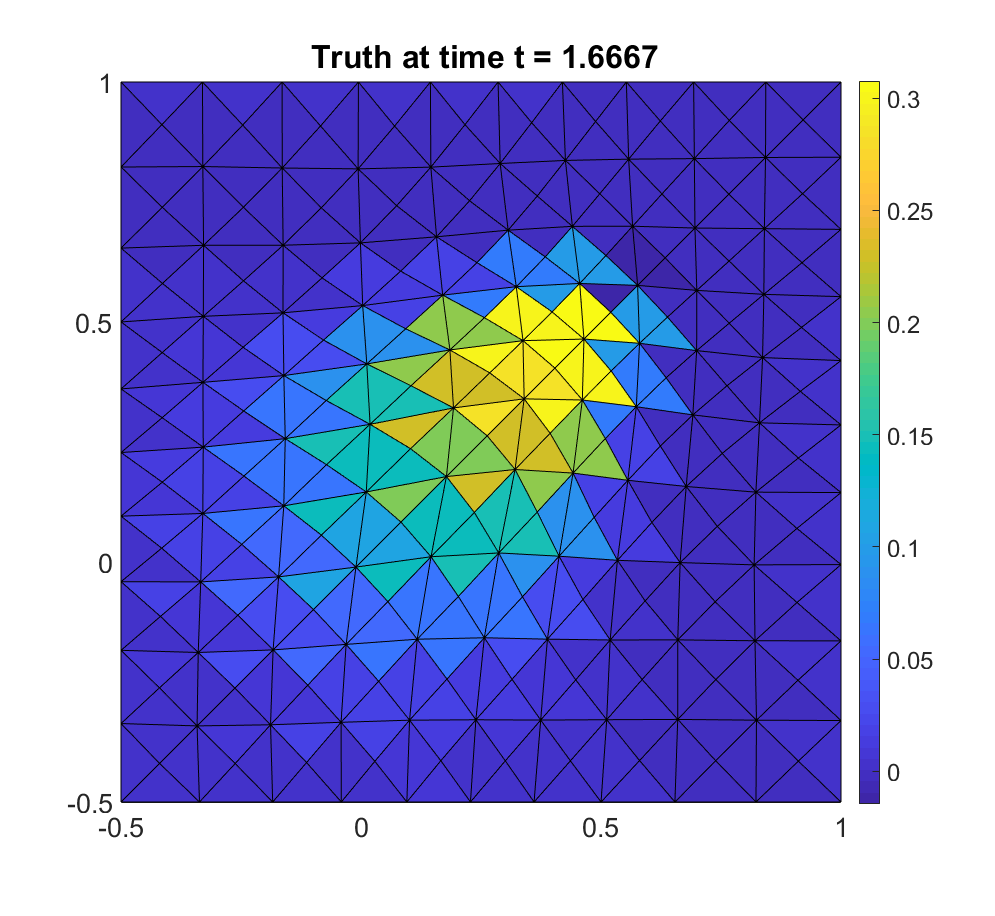}\hfill
     \includegraphics[width=.24\textwidth]{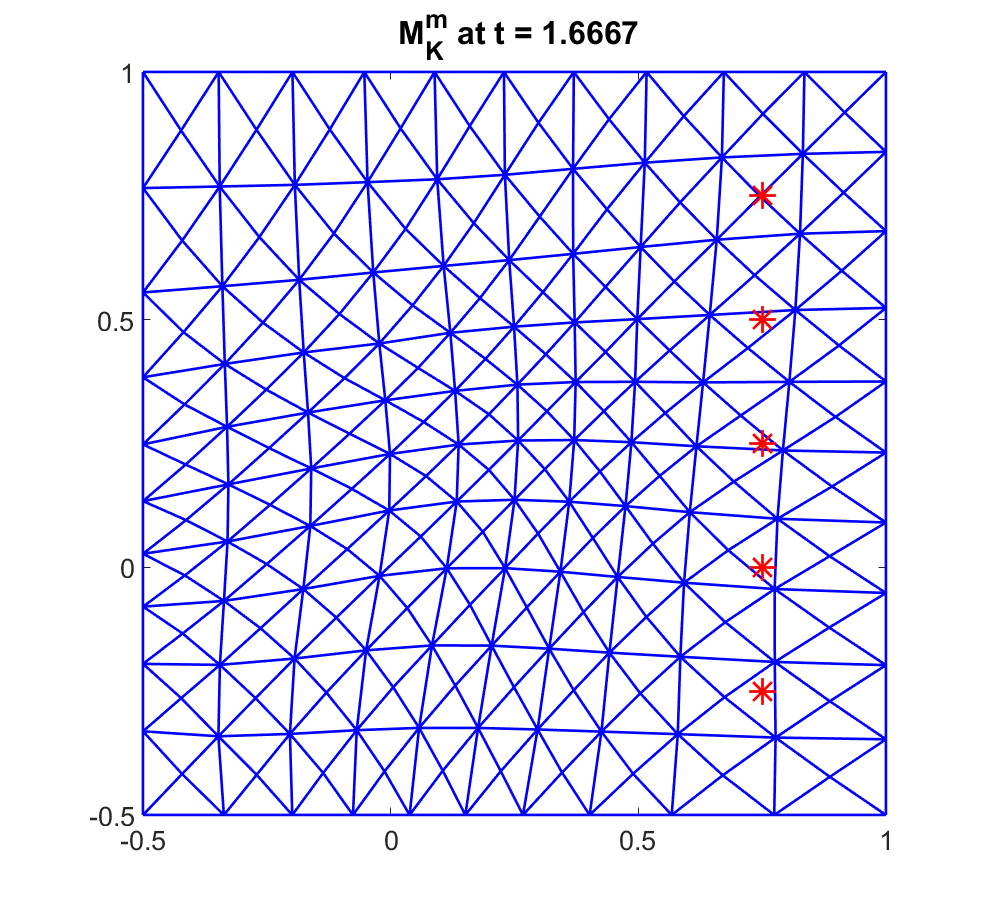}\hfill
     \includegraphics[width=.24\textwidth]{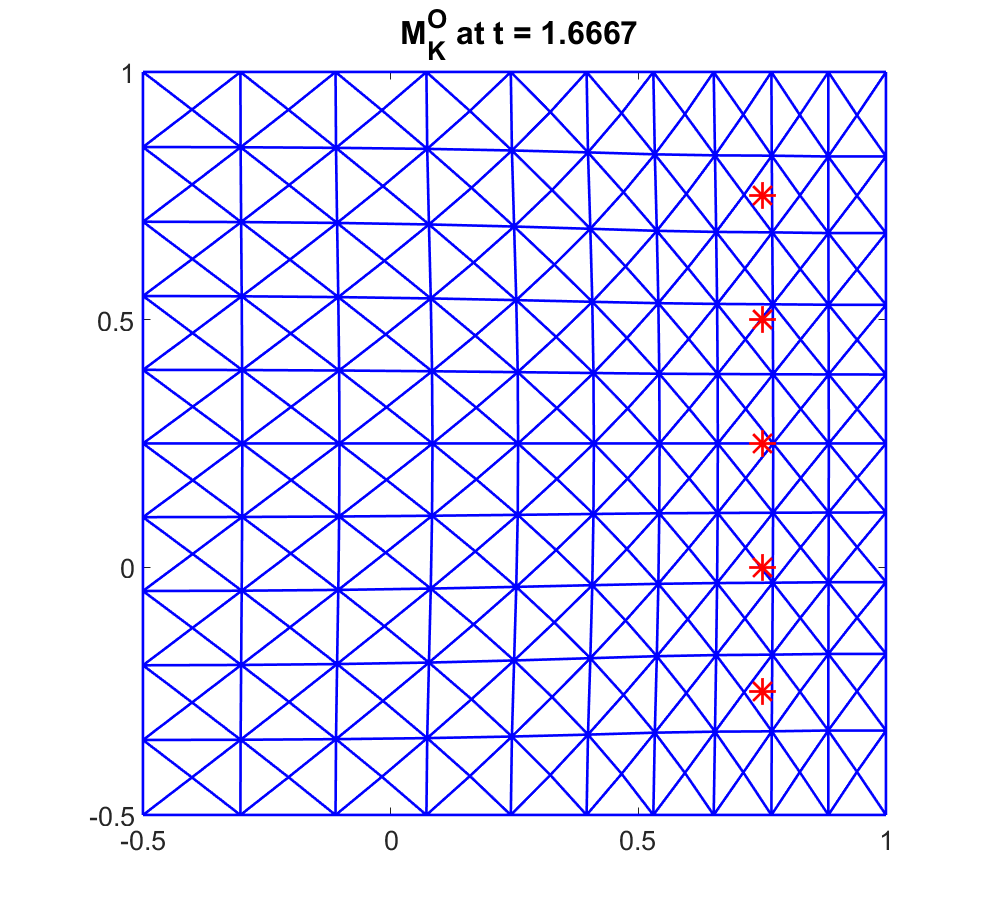}\hfill
     \includegraphics[width=.24\textwidth]{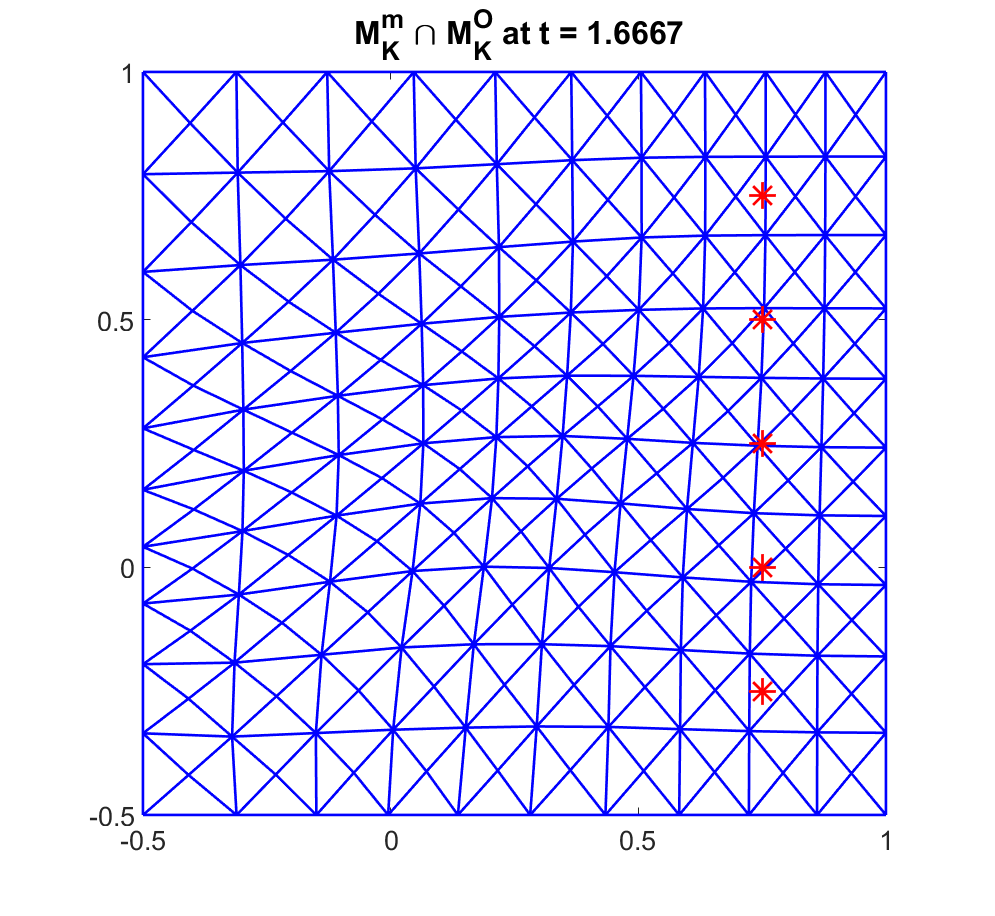}
     \caption{A comparison of the different choices of common mesh.}
     \label{fig:2dmeshes}
 \end{figure}

 The choice of concentration with the common mesh is easily seen in 2 dimensions. Figure \ref{fig:2dmeshes} shows a snapshot of the numerical solution of the 2D inviscid Burgers equation \eqref{eq:Burgers2d} at time $t = 1.6667$. At this time, the front has fully formed and is beginning to propagate to the upper left corner. There are five observations sampled along the line $x=0.75$. The truth is shown in the leftmost panel, followed by the three choices for the common meshes computed from $\M_K^m$, $\M_K^O$, and $\M_K^m \cap \M_K^O$. The common mesh taken from $\M_K^m$ focuses the mesh points near the front, while the common mesh from $\M_K^O$ concentrates the mesh near the line of observations. The choice on the far right ($\M_K^m \cap \M_K^O$) concentrates the mesh near both of these attributes.

 \subsubsection{Experiment 5: Tuning localization and inflation}
\begin{figure}[ht]
    \centering
    \includegraphics[width = 0.3\textwidth]{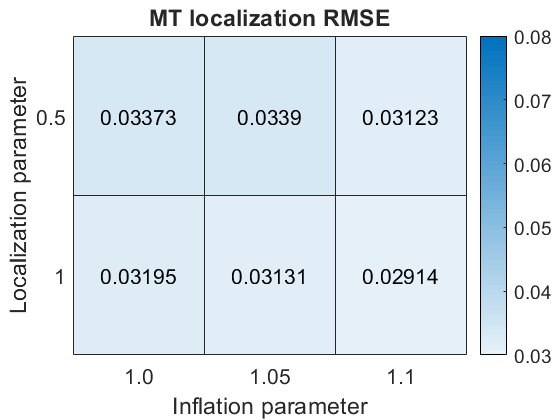}\hfill
    \includegraphics[width=0.3\textwidth]{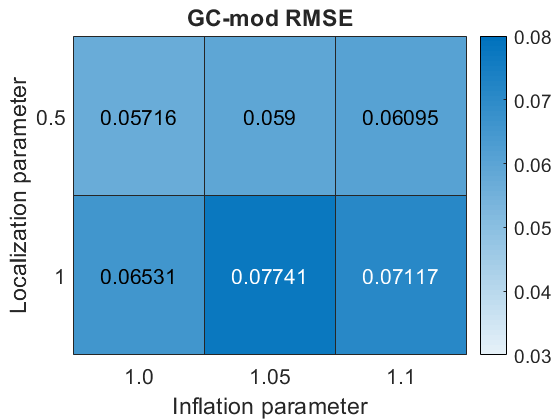}\hfill
    \includegraphics[width=0.3\textwidth]{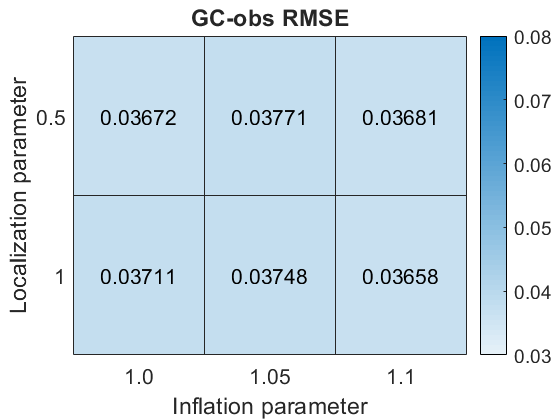}
    \caption{Results for simultaneous tuning of localization and inflation for the MT, GC-mod, and GC-obs localization schemes.}
    \label{fig:tune_rho_irad2D}
\end{figure}

 Like in the 1-dimensional case, the localization parameters for the various localization schemes are tuned simultaneously with the inflation parameter. Here we consider a reduced parameter space, using an inflation factor $\rho \in \{1,1.05,1.1\}$ and localization parameter $L\in\{0.5,1.0\}$. The results are presented in Figure \ref{fig:tune_rho_irad2D}.
Just like in the 1D case, the MT localization scheme is robust with respect to the tuning of the localization parameter.

 \subsubsection{Experiment 6: Localization scheme comparison}
 Consider the following localization schemes, tuned in Experiment 5:
\begin{itemize}
\item MT localization with multiplicative inflation parameter $1.1$ and localization parameter $1.0$.
\item GC-mod localization with multiplicative inflation parameter $1.0$ and localization parameter $0.5$.
\item GC-obs localization with multiplicative inflation parameter $1.1$ and localization parameter $1$.
\end{itemize}

\begin{figure}[ht]
    \centering
    \includegraphics[width=.7\textwidth]{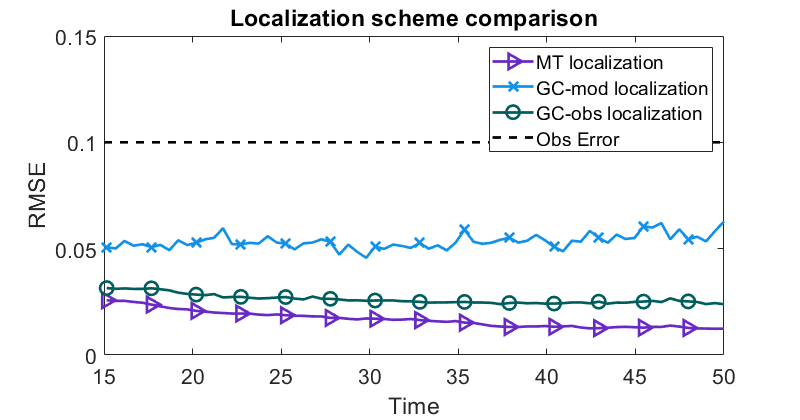}
    \caption{Time series RMSE comparison for the different localization schemes for 2D inviscid Burgers. Each trial uses the optimal localization parameters that were tuned in Experiment 5.}
    \label{fig:locRMSE2D}
\end{figure}

The time series RMSE of each of these localization schemes in shown in Figure \ref{fig:locRMSE2D}. In these experiments, $\Sigma = 0.01I$, so the localization scheme is considered a success if it is on the order of $0.1$. 

 \subsubsection{Experiment 7: Choice of common mesh}
\begin{figure}[h]
    \centering
    \includegraphics[width = .8\textwidth]{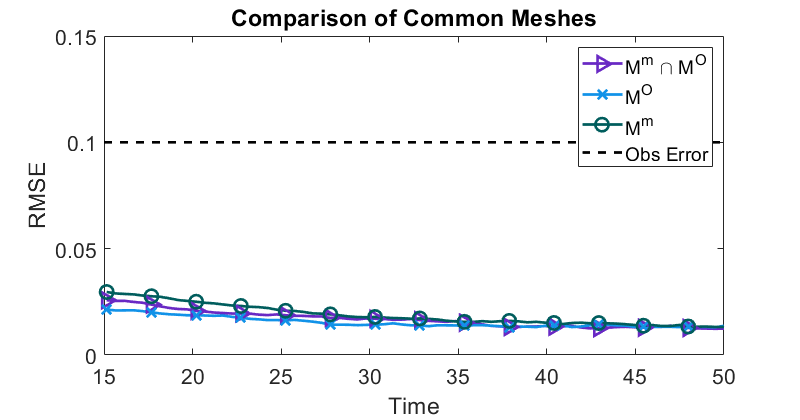}
    \caption{Time series RMSE for tuned DA with 2D inviscid Burgers equation using different choices for the common mesh.}
    \label{fig:mesh_rmse2D}
\end{figure}

 Again we compare the choice of the common mesh for the tuned MT localization scheme: $\M \in \{\M^{O},\M^{m},\M^{m}\cap\M^{O}\}$. A comparison of these common meshes is shown in Figure \ref{fig:mesh_rmse2D}, using the tuned localization and inflation parameters found in Experiment 5. All choices for the common mesh yield a stable DA procedure. 

\subsubsection{Experiment 8: Comparison of error covariances}
\begin{figure}[h]
    \centering
    \includegraphics[width = .85\textwidth]{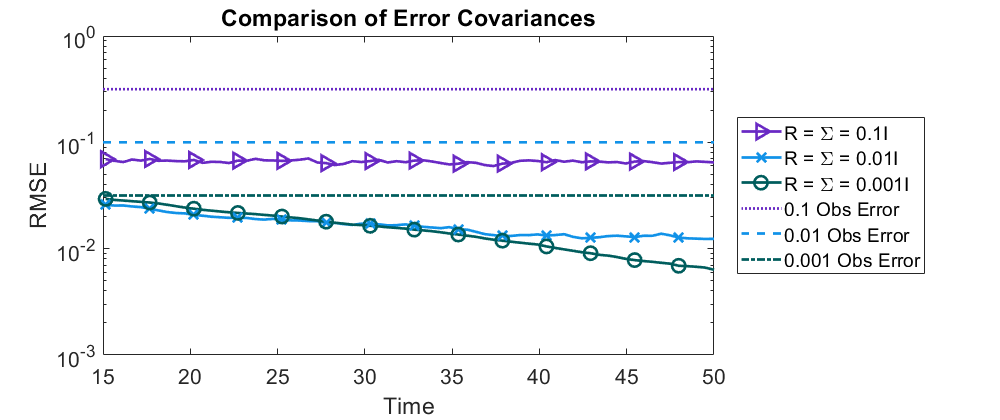}
    \caption{Time series RMSE for 2D Burgers with tuned MT localization and different error covariances.}
    \label{fig:covariances2D}
\end{figure}

Using the tuned MT localization scheme from Experiment 5, we test the robustness of this DA procedure on the 2D inviscid Burgers \eqref{eq:Burgers2d} with different error covariances. For each experiment, we set the model error, observation error, and initial error covariances equal to a scalar multiple of the identity.  The results in Figure \ref{fig:covariances2D} show stable RMSE. 

\subsubsection{Experiment 9: 2D Burgers with Noisy Data}
\begin{figure}[h]
    \centering
    \includegraphics[width = .85\textwidth]{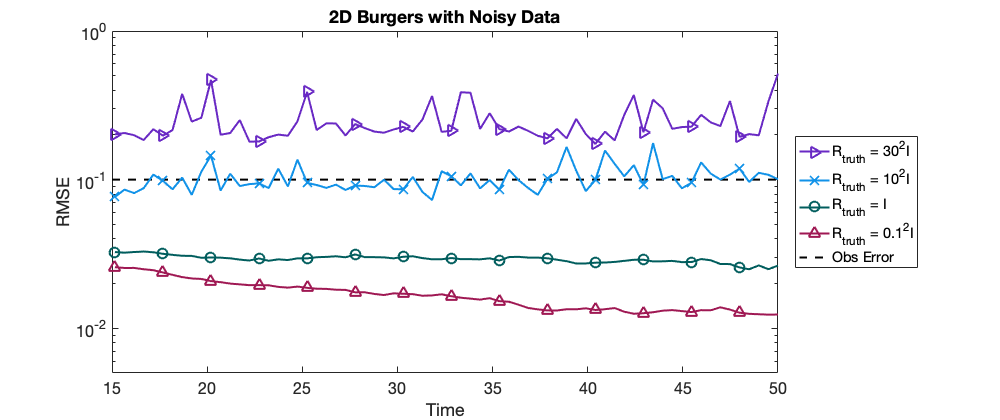}
    \caption{Time series RMSE for 2D Burgers with noisy data. It is assumed that the observation error has a normal distribution with mean $0$ and covariance matrix assumed to be $R=0.01I$, but the synthetic observations are produced from the truth using an error covariance matrix $R_{truth}$.}
    \label{fig:baddata}
\end{figure}

To test this configuration in a more difficult regime, consider \eqref{eq:Burgers2d} with noisy data. That is, suppose the data error covariance is believed to be $R = 0.01 I$, but the data is actually sampled with an error covariance of $R_{truth} = \alpha^2 I$, where $\alpha \geq 0.1$. For large values of $\alpha$, (e.g., $\alpha\geq 100$), the analysis update causes discontinuities in the numerical solution that the time integrator cannot overcome. For moderately large values of $\alpha$ (e.g., $\alpha = 50$) this is sometimes an issue, but often is not. We have found that our approach works consistently for $\alpha \leq 30$, and produces a stable DA algorithm for these values. Note that in these experiments, the initial perturbation and the model error continue to be sampled from $\Sigma = 0.01I$. The results for an observational time step of $\Delta t = 0.5$ are shown in Figure \ref{fig:baddata}. 

\subsubsection{Experiment 10: Choice of Interpolation}
\begin{figure}
    \centering
    \includegraphics[width=0.7\textwidth]{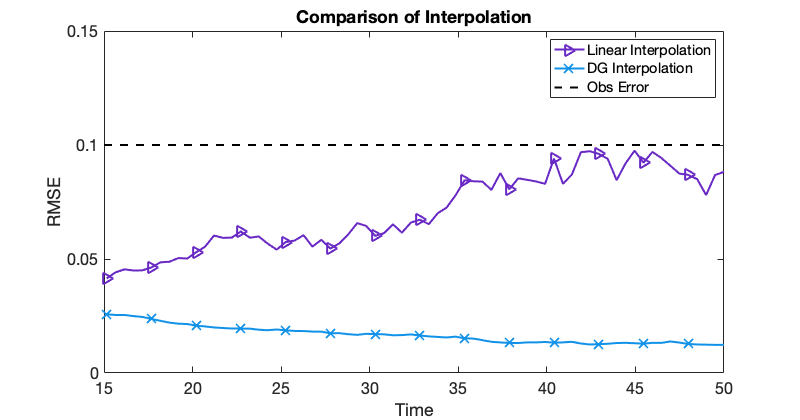}
    \caption{Linear interpolation vs. DG interpolation for DG discretization of 2D inviscid Burgers.}
    \label{fig:interp}
\end{figure}

In all of the experiments above, we use the DG interpolation, as detailed in Section \ref{sec:details-dg}. There is not much difference in accuracy in the 1D case, but the choice of interpolation makes a significant difference in the 2D case. Figure \ref{fig:interp} shows the difference in RMSE when using the linear interpolation instead of the DG interpolation.

\section{Discussion and Conclusions}
Through the use of an adaptive common mesh, we develop an ensemble based DA scheme where each of the ensemble members evolve independently on their own adaptive meshes. At each observational timestep, the ensemble members are interpolated to the adaptive common mesh, updated according to the DA scheme, and then interpolated back to their individual meshes. 

We follow the MMPDE adaptive meshing strategy where the mesh of each ensemble member is determined by a matrix-valued monitor function, also called a metric tensor, so that the mesh is viewed as uniform in that metric. At each observational timestep, an adaptive common mesh is calculated. There are several choices for this common mesh. One choice, $M^m$, is obtained by taking the intersection of the ensemble members' metric tensors. This results in a common mesh that in some sense satisfies all of the ensemble members. Another option is to concentrate the common mesh near observation locations or observation trajectories. Concentrating the mesh near the observation locations reduces the amount of interpolation error at each observational timestep. 
A third choice is to intersect $\M^m$ with $\M^O$. Using the observational mesh $\M^O$ in a DA scheme reduces the transient time in converging to the asymptotic behavior, but regardless of which is employed as the common mesh, all choices $\M^m$, $\M^O$, and $\M^m\cap\M^O$ produce stable results. The efficacy of several techniques developed in this work is illustrated using sharp interface problems, in particular 1D and 2D inviscid Burgers equations, under a discontinuous Galerkin discretization.

We develop a new adaptive localization algorithm based on the metric tensor of the common mesh $\M^m$. The MT adaptive localization scheme uses the metric tensor to define a domain localization strategy that is dynamically updated in time and space. For the 1D and 2D inviscid Burgers equations, the MT localization scheme compared favorably with the GC localization schemes. One of the benefits of the MT localization scheme is that it is robust with respect to the tuning parameters and requires less precise tuning than GC localization in either the model space or in the observation space.

The interpolation that is used at each observational timestep can have a significant impact on the performance of the DA scheme. Using a DG discretization for the PDE together with a DG-based interpolant allows the ensemble members to maintain the advantages of DG discretization independent of their supporting mesh.  For the 1D inviscid Burgers problem, there was no significant difference in RMSE between linear and DG interpolation. For the 2D inviscid Burgers equation, however, using a DG-based interpolation scheme improved the overall performance of the DA scheme as compared to linear interpolation.

This metric tensor approach to data assimilation on adaptive moving meshes, as well as the MT localization scheme, is applicable in higher spatial dimensions.  There are several interesting avenues for further investigation. These include the development of adaptive meshes in which a single (average) mesh supports all ensemble members over an observation cycle, further development of adaptive meshing techniques to minimize error due to uncertainties in the location of observations, and the development of goal-oriented meshing functionals specifically designed to increase the skill of the data assimilation scheme.


\newpage
\appendix
\section{Notation Glossary}
\begin{tabular}{c|l}
Variable & Description\\
\hline
     $\mathcal{C}$ & Gaspari-Cohn localization function\\
     $D$ & dimension of observation space\\
     $d$ & spatial dimension of PDE\\
     $e_i$ & $i^{th}$ ensemble member\\
     $F_K$ & affine mapping between elements in mesh and reference element\\
     $G$ & Mesh function\\
     $\mathcal{H}$ & (nonlinear) observation operator\\
     $H$ & linearization of $\mathcal{H}$\\
     $h$ & refers to spatial discretization, as a subscript\\
     $I$ & identity\\
     $j$ & spatial discretization index \\
     $\mathbf{K}$ & Kalman gain matrix\\
     $K$ & element in mesh\\
     $\hat K$ & reference element\\
     $\mathcal{L}$ & function whose gradient gives mesh equation\\
     $L$ & localization parameter\\
     $M$ & dimension of state space \\
     $\hat m$ & predicted ensemble mean\\
     $\M$ & metric tensor\\
     $\M_K$ & metric tensor at element $K$\\
     $\M_K^{e_i}$ & metric tensor of $i^{th}$ ensemble member at element $K$\\
     $\M_K^m$ & metric tensor of mesh obtained from metric tensor intersection of ensemble meshes at element $K$\\
     $\M_K^O$ & metric tensor of observation mesh at element $K$\\
     $n$ & time index\\
     $\mathcal{N}$ & normal distribution\\
     $N_e$ & number of ensemble members \\
     $N_O$ & number of observations\\
     $N$ & number of simplicial complexes in mesh\\
     $P^b$ & background error covariance matrix \\
     $R$ & observation error covariance matrix \\
     $S$ & length of spatial domain for 1D Burgers\\
                    \end{tabular}

     \begin{tabular}{c|l}
     Variable & Description (continued)\\
     \hline
     $t$ & time variable\\
     $\Delta t$ & observation time scale\\
     $T$ & matrix transpose, as a superscript\\
     $T$ & time index at final time\\
     $\mathcal{T}_h^{e_i}$ & $i^{th}$ ensemble mesh\\
     $\mathcal{T}_h$ & common mesh\\
     $u$ & state variable\\
     $\hat u$ & state forecast \\
     $u^b$ & background state \\
     $U$ & numerical solution \\
     $V$ & state vector augmented with mesh locations \\
     $x$ & node location in 1D\\
     $x^O$ & observation location\\
     $y$ & observation variable\\
     $z$ & continuous time DA observation variable \\
     $\Psi$ & physical model for discrete time dynamical system\\
     $\xi$ & model noise \\
     $ \eta$ & observation noise \\
     $\Sigma $ & model error covariance matrix\\
     $\tau$ & scalar to determine mesh velocity \\
     $\Omega$ & polyhedral domain\\
\end{tabular}

    \bibliographystyle{abbrv}
    \bibliography{KU_MMDA}

\end{document}